\RequirePackage{fix-cm}

\documentclass{svjour3}                     
\smartqed  

\usepackage{color}
\usepackage{verbatim}
\usepackage{hyperref}
\usepackage{booktabs}
\usepackage{siunitx}
\usepackage{fancyvrb}
\usepackage{graphicx}
\usepackage{comment}
\usepackage[ruled,vlined]{algorithm2e}
\usepackage{mathtools}
\usepackage{graphicx}
\usepackage{subcaption}
\usepackage{multirow}
\usepackage{caption}
\usepackage{cite}
\usepackage{booktabs, dcolumn}

\newcommand{\norm}[1]{\left\lVert#1\right\rVert}

\newcommand{\x}{\mathbf{x}}
\newcommand{\y}{\mathbf{y}}
\newcommand{\zz}{\mathbf{z}}
\newcommand{\A}{\mathbf{A}}

\newcommand{\D}{\mathbf{D}}
\newcommand{\X}{\mathbf{X}}
\newcommand{\bb}{\mathbf{b}}
\newcommand{\dd}{\mathbf{d}}
\newcommand{\thh}{\pmb{\theta}}
\newcommand{\Gam}{\pmb{\Gamma}}

\usepackage{mathptmx}      
\begin{document}

\title{Distributed Primal Outer Approximation Algorithm for Sparse Convex Programming with Separable Structures
}

\titlerunning{Distributed Primal Outer Approximation Algorithm for Sparse Convex Programming}        

\author{Alireza Olama          \and       
        Eduardo Camponogara     \and   
        Paulo R. C. Mendes
}


\institute{Alireza Olama \at
              Department of Automation and System Engineering, Federal University of Santa Catarina, Florian\'{o}polis, Brazil \\
              \email{alireza.olama@posgrad.ufsc.br}           
           \and
           Eduardo Camponogara \at
             Department of Automation and System Engineering, Federal University of Santa Catarina, Florian\'{o}polis, Brazil \\
              \email{eduardo.camponogara@ufsc.br}  
              \and
           Paulo R. C. Mendes \at
             Fraunhofer Institute  for Industrial Mathematics, Kaiserslautern, Germany \\
              \email{paulo.mendes@itwm.fraunhofer.de}  
}

\date{Received: date / Accepted: date}

\maketitle

\begin{abstract}
This paper presents the \textbf{D}istributed \textbf{P}rimal \textbf{O}uter \textbf{A}pproximation (DiPOA) algorithm for solving Sparse Convex Programming (SCP) problems with separable structures, efficiently, and in a decentralized manner. The DiPOA algorithm development consists of embedding the recently proposed Relaxed  Hybrid Alternating Direction Method of Multipliers (RH-ADMM) algorithm into the Outer Approximation (OA) algorithm. We also propose two main improvements to control the quality and the number of cutting planes that approximate nonlinear functions. In particular, the RH-ADMM algorithm acts as a distributed numerical engine inside the DiPOA algorithm. DiPOA takes advantage of the multi-core architecture of modern processors to speed up optimization algorithms. The proposed distributed algorithm makes practical the solution of SCP in learning and control problems from the application side. 
This paper concludes with a performance analysis of DiPOA for the distributed sparse logistic regression and quadratically constrained optimization problems.  Finally, the paper concludes with a numerical comparison with state-of-the-art optimization solvers.
\subclass{MSC code1 \and MSC code2 \and more}
\end{abstract}

\section{Introduction}
\begin{table}[]
    \caption{Acronyms}
    \centering
    {    \begin{tabular}{ll}
    \toprule
       Acronym & Meaning\\
        \toprule
       CN &Communication Network\\
       DiPOA& Distributed Primal Outer Approximation\\
       D-MILP & Distributed Mixed Integer Linear Program\\
       D-NLP & Distributed Nonlinear Program \\
       D-MINLP & Distributed MINLP\\
              DSLR & Distributed Sparse Logistic Regression\\
       ET-SoCut& Event Triggered SoCut\\
       GBD & Generalized Benders Decomposition \\
        LFC & Local Fusion Center\\
        MIP & Mixed Integer Programming \\
       MPI & Message Passing Interface\\
       MINLP & Mixed Integer Nonlinear Program\\
       MILP & Mixed Integer Linear Program\\
       MIQP & Mixed Integer Quadratic Program\\
       MIQCP & Mixed Integer Quadratically Constrained Program\\
       NLP & Nonlinear Programming\\
        OA & Outer Approximation\\
       RH-ADMM & Relaxed-Hybrid Alternating Direction Method of Multipliers\\
       SCP & Sparse Convex Programming\\    
        SLR & Sparse Logistic Regression\\
       SoCut & Second Order Cut \\
       SQCQP& Sparse Quadratically Constrained Quadratic Programming\\

    \bottomrule
    
    \end{tabular}
    }
    {}
\end{table}
A broad class of modern real-world applications of engineering and artificial intelligence problems consists of mathematical optimization problems with a constraint that allows only up to a certain number of variables to be nonzero.
This form of constraint is called a \textit{sparsity constraint} (or a cardinality constraint) and any convex optimization problem containing such a constraint is referred to as a \textit{Sparse Convex Programming} (SCP) problem. 
Studied by \cite{Bienstock1996}, one of the first  SCP problems involved a quadratic objective function and a set of linear constraints.
The SCP problems arise in various applications such as linear and nonlinear compressed sensing \cite{blumensath2009iterative,blumensath2013compressed}, best subset selection on regression and classification problems \cite{Bertsimas2015,Bertsimas2017,Bertsimas2017b,bertsimas2021sparse}, mixed-integer programs \cite{wang2020new,bourguignon2015exact}, and so on. See excellent survey for more details about applications \cite{tillmann2021cardinality}.
Although SCP problems find uses in a variety of applications, they are hard to solve for global optimality. Since the sparsity constraint is a union of finitely many sub-spaces, the SCP problem is a combinatorial optimization and hence is generally NP-hard \cite{Wang2020,Sun2013,Bai2016}.

Due to the numerical difficulties imposed by the sparsity constraint, most of the SCP solution algorithms try to reformulate or approximate the original problem such that it becomes simpler to handle.
Accordingly, most of the algorithmic solutions for SCP are often divided into two categories, namely, approximation and exact methods. The approximation methods often try to provide a good convex representation of the sparsity constraint and obtain the desired sparsity by solving a sequence of convex optimization problems. One of the most popular approximation methods is the $\ell_1$ convex relaxation method which provides a convex lower bound on the sparsity constraint in the unit $\ell_\infty$ norm. Because of their efficiency, the approximation methods are popular in solving feature selection problems found in the statistical learning literature \cite{tibshirani1996regression,boyd2011distributed}. There exist also certain conditions (\textit{e.g.}, restricted isometry property condition) that suffice for the approximation methods to find the exact solution \cite{candes2008restricted}.

The exact methods, on the other hand, try to view the SCP problem as a Mixed Integer Nonlinear Programming (MINLP) problem since it can be reformulated in such a way by introducing suitable binary variables \cite{Bertsimas2015,Bienstock1996}. In this paper, we mainly focus on the exact methods.
As the first attempt, \cite{Bienstock1996} proposed a Mixed Integer Quadratic Programming (MIQP) problem to reformulate the SCP problem with a quadratic objective and linear constraints. The problem is then solved by employing a tailored branch and cut algorithm.
\cite{Bertsimas2009} proposed a tailored algorithm based on the branch and bound method to solve the SCP problem appearing in sparse linear regression and portfolio selection problems. 
\cite{Aguilera2017} provided a MIQP formulation for solving a sparsity constrained model predictive control problem. Although the methodology works efficiently for the given application, it lacks generalization to more realistic scenarios.
As another instance, \cite{bertsimas2022scalable} considered a scalable algorithm to solve cardinality constrained portfolio optimization problems. This paper also proposed a MIQP framework along with multiple improvements and heuristics such as high-quality warm-starts, a prepossessing technique, and so on.
\cite{aytug2015feature} proposed a MIQP based framework to solve a sparsity constrained support vector machine problem. In this paper, the authors utilized the Generalized Benders Decomposition (GBD) algorithm to select the best subset of features during the model training.
As another application, \cite{Bertsimas2015,Bertsimas2009} proposed an efficient algorithm to obtain an exact solution for the best subset selection problem in linear regression. In this paper, the authors proposed a MIQP reformulation for the resulting SCP problem. The resulting MIQP problem is then solved by a tailored branch and bound algorithm. Although the MIQP approach presented in \cite{Bertsimas2015,Bertsimas2009} obtains an exact solution to the best subset selection problem, however, there are cases where the MIP approach performs worse than Lasso both in terms of solution quality and, in computational time as it is reported in \cite{hastie2020best}.
As another example, an MINLP model to solve sparse classification problems is proposed by \cite{Bertsimas2017,Bertsimas2017a}. The resulting MINLP model is then solved by utilizing the Outer Approximation (OA) algorithm \cite{Grossmann2002,Kronqvist2019}.
 \cite{kamiya2019feature} consider a feature subset
selection problem for the multinomial logistic model with $L_2$ regularization. This problem is then transformed into a SCP problem which in turn is solved with the OA algorithm.
Additionally, \cite{bertsimas2021sparse} provides a MINLP formulation to solve a sparse convex regression problem.
Besides the decomposition methods (\textit{e.g.}, OA algorithm), some works focused on single-tree methods such as branch and bound, branch and cut, and branch and price algorithms.
For example, \cite{wang2020new} used the branch and bound method, with domain cut and partition scheme, to solve a quadratic sparsity constrained problem that has application to portfolio optimization. 
As another approach, \cite{sant2020solving} proposed an MINLP modeling framework along with a SCP problem to model and solve the index tracking problem. The problem is then solved by a variant of the branch and cut algorithm. 

There exist relatively few works where a general convex MINLP model is used to reformulate the SCP problem. For instance, instead of formulating the SCP problem for a particular application, \cite{Bai2016} considered a general convex formulation for the SCP problems. The authors proposed a splitting augmented Lagrangian method to solve the resulting SCP problem. The efficiency of the proposed algorithm is then tested on portfolio selection and compressed sensing problems.

\subsection{Limitations of the Existing Works}
As mentioned, most of the works focused on particular cases and applications of the general SCP problem whose objective function is either linear or quadratic while enforcing linear constraints. Although there has been considerable progress in solving the SCP problems with quadratic and linear objectives, in many real-world applications the SCP problems often consist of a convex nonlinear objective and constraints. Up to the authors' knowledge there exist a few works that consider a general convex MINLP (\textit{e.g.}, \cite{Bai2016}). Additionally, most of the real-world applications of the SCP are inherently distributed or can be easily decomposed into smaller problem instances which can be solved distributedly \cite{nedic2018distributed}.
As an example, in some instances of regression and classification problems, the dataset is distributed over a network and the data privacy policy does not allow the collection of the entire dataset in a single point.
There are similar situations in various applications, such as energy management of renewable power systems, distributed compressed sensing, and so on. 
Hence, a distributed solution is inevitable. It is worth also noting that although the distributed reformulation of the SCP problem is necessary, it imposes more complexity on the problem due to the coupling behavior in either objective function or constraints.
Accordingly, we consider the limitation of the existing works as,
\begin{enumerate}
    \item a general convex formulation to handle nonlinearities.
    \item a distributed formulation to handle coupling behaviors.
\end{enumerate}
\textcolor{black}{The main motivation of this paper is to overcome the limitations of the existing works and provide an algorithmic framework to solve the SCP problem, efficiently and distributedly.}

\subsection{Main Contributions}

Despite the advantages of the centralized architecture, it might not be the best solution strategy in some real-world applications where the problem is inherently distributed. Some important instances are found in estimation, decision making, learning, and control applications \cite{notarstefano2019distributed,nedic2018distributed,nedic2015distributed}. Essentially, most of these modern problems are defined over a so-called Communication Network (CN)  in which agents aim at cooperatively solving complex problems by local computation and communication \cite{notarstefano2019distributed}. To solve optimization problems over CNs, it might not be trivial to apply centralized optimization algorithms, which require the data to be managed by a single entity. In such a situation, the problem data is often distributed over the network, and it is usually undesirable to collect them at a unique node. To this end, \textit{parallel computing} serves as a source of inspiration. In the context of convex MINLP problems, there exist a few works in the literature that present a parallelizable version of the OA algorithm. As an example \cite{muts2020decomposition} presents a two-phase algorithm for solving convex MINLP problems, called decomposition-based outer approximation algorithm. In the first phase, a sequence of linear integer relaxed sub-problems is solved while in the second phase, the algorithm solves a sequence of mixed-integer linear programming sub-problems. However, To the best of our knowledge, there are no approaches in the literature that consider a distributed set-up for the SCP problems over CNs.

This paper proposes an algorithm, the so-called \textit{Distributed Primal Outer Approximation (DiPOA)}, to distributedly solve SCP problems over a CN. In principle, DiPOA extends the OA algorithm by embedding a fully decentralized algorithm, namely the Relaxed Hybrid Alternating Direction Method of Multipliers (RH-ADMM) which was proposed by the authors in \cite{olama2019}. The RH-ADMM assumes a special architecture on the CN, the so-called \textit{hybrid architecture}, which was developed to solve distributed convex optimization problems. Additionally, to further improve the efficiency of DiPOA, multiple heuristics and improvements are proposed. The main contributions of this work are:
\begin{enumerate}

\item Equivalent distributed MINLP models for large-scale SCP problems. This is achieved by utilizing the hybrid architecture for the CN architecture and consensus optimization concepts.

\item DiPOA, an MINLP algorithm to solve the SCP problems, which is designed by embedding the RH-ADMM into the OA algorithm and performing distributed computations.

\item A specialized distributed feasibility pump method to accelerate the convergence of the DiPOA.

\item A practical distributed method for the RH-ADMM algorithm to detect problem infeasibility before performing the main computational steps.

\item A distributed second-order cut and an event-triggered cut generation scheme to improve performance and computational efficiency.

\item A performance analysis of DiPOA regarding its applications to solve Distributed Sparse Logistic Regression (DSLR) and Sparse Quadratically Constrained Quadratic Programming (SQCQP) problems over a CN network. 


\end{enumerate}
\subsection{Paper Organization}
This paper is organized as follows; Section \ref{problem-back} consists of problem formulation and backgrounds. The solution architecture and the MINLP formulations are provided in  \ref{reformlulation}. Section \ref{dipoa} presents the development of the DiPOA algorithm, algorithm improvements, and heuristics. The implementation details and the numerical experiments are explained in Section \ref{implementation} and \ref{exper}, respectively. Finally, Section \ref{concolusions} draws some conclusions and suggests research directions for future works.
\section{Problem Formulation and Background}\label{problem-back}

This section presents the formulation of the considered SCP problem and explains core concepts of the RH-ADMM algorithm.
\subsection{Problem Definition}
This paper aims to study the SCP problems with separable structures in the following form,
\begin{subequations}\label{dis-ccp}
\begin{align}
\text{P$_1$: }\qquad    \min_{\x \in \mathcal{R}^n}\,\,             &\sum_{i=1}^{N}f_i(\x) \\
  s.t.~  &g_h(\x) \leq 0, \,\, \forall h=1,...,m \\
    & \x \in \Omega \\
    &\norm{\x}_0 \leq \kappa \label{D-CCP}
\end{align}
\end{subequations}
where $\x\in \mathcal{R}^n$ is a global vector of decision variables which is shared between each node, $N$ is the number of the nodes of the CN, $f_i: \mathcal{R}^n\xrightarrow{}\mathcal{R}$ and $g_h : \mathcal{R}^n\xrightarrow{}\mathcal{R}$ are continuously differentiable functions that describe the local objective functions and nonlinear constraints, respectively. Note that problem \eqref{D-CCP} is separable with respect to the local data that constructs each objective function, however, the existence  of a global decision vector $\x$ provides a coupling between local objective functions.
 
 The sparsity constraint, which counts the number of non-zeros of the vector $\x$, is modeled by the $\ell_0$ norm, being defined as $\norm{(x_1,...,x_n)}_0 = \sum_{i=1}^n 1 (x_i \neq 0)$. A scalar $\kappa < n$ defines the maximum number of non-zero elements of the decision vector $\x$. Moreover, $\Omega = \{\mathbf{x}\in \mathcal{R}^n\ |\ \D \x \leq \dd,\, \A\x=\bb\}$ is a polytope which consists of local linear equality and inequality constraints of the problem. As detailed in Remark \ref{remark_1}, $P_1$ is distributed in the sense that each function $f_i(\cdot)$ is computed only by node $i$ of the network.
\\
\textbf{Assumption: }
The local objectives $f_i, i\in\{1,..., N\},$ and constraints $g_h, h\in\{1,...,m\},$ are assumed to be convex, closed, and proper functions.
\begin{remark} \label{remark_1}
There exist several distributed settings for optimization problems. In one of the most common formulations, the decision variables are shared between the nodes, and the problem data provided for each node of the CN is considered to be private. For this reason,  as problem \eqref{dis-ccp} indicates, each node has only enough information to construct its objective function, while the decision variable vector is shared among all nodes. This distributed setting is a common pattern in some significant learning and control applications.
For instance, in distributed ML, an efficient technique to deal with large volumes of data consists of distributing the data over a network. In this case, a significant reduction in memory consumption can be achieved while keeping the same unknown parameters of the model.
\end{remark}

\textcolor{black}{In the context of mathematical optimization, the term consensus refers to a class of problems with $N$ separable local objective functions $f_i: \mathcal{R}^n\xrightarrow{} \mathcal{R}$ and vectors of decision variables $\x_i \in \mathcal{R}^n$ subject to a set of equality constraints, the so-called \textit{consensus} or \textit{consistency} constraints \cite{boyd2011distributed, notarstefano2019distributed}. The consensus constraints enforce the local decision variables to converge to a global variable $\y \in \mathcal{R}^n$ and usually such constraints are defined as $\x_i = \y$, for $i = 1,\dots, N$. The solution algorithms usually impose the consensus constraints iteratively until the constraints $\x_i = \y$, for $i = 1,\dots, N$ are satisfied within a given tolerance.
}

\textcolor{black}{The main goal of this work is to develop an efficient numerical optimization algorithm for distributedly solving the SCP problem \eqref{dis-ccp}.} To do so, we will first transform the SCP problem \eqref{dis-ccp} into a consensus optimization problem\footnote{See problem (\ref{admm-ccp}a)-(\ref{admm-ccp}e).}  by imposing the hybrid architecture to the CN \cite{ma2018fast,olama2019}. This is achieved by considering multiple nodes of the CN as coordination nodes. Such a hybrid architecture allows us to distribute the problem data over multiple nodes when a big data set is considered. Moreover, when the problem is inherently distributed and its data is already distributed over a CN, the hybrid architecture provides a framework to perform the computations locally without using a central coordinator. 

Following the hybrid architecture for CNs, we will reformulate the consensus SCP problem as a distributed MINLP which is equivalent to the original SCP problem \eqref{dis-ccp}. The complexity of the distributed MINLP is attributed to its binary variables.
\subsection{A Short Summary of MINLP Methods}
Generally speaking, most of the modern MINLP algorithms and methods apply a decomposition scheme in which the original MINLP problem is decomposed into two sub-problems, namely, \textit{master's} and \textit{primal} problems.
\textcolor{black}{Generally, the master's problem is a linear (mixed) integer problem that is updated dynamically and the primal problem is a nonlinear problem with continuous decision variables.}
\textcolor{black}{MINLP algorithms such as Generalized Benders Decomposition (GBD) \cite{geoffrion1972generalized}, Outer Approximation (OA) \cite{duran1986outer}, Extended Cutting Plane (ECP) \cite{westerlund1995extended}, and Extended Supporting Hyperplane (ESH) \cite{kronqvist2016extended} follow the same strategy. } For a detailed review of MINLP algorithms and methods, refer to the survey by \cite{Kronqvist2019}.
%

Despite the advantages mentioned above, duality-based methods such as GBD can perform poorly in certain situations.

\begin{remark}
Despite several algorithms for solving MINLP problems being available off the shelf, few methodologies solve such challenging problems in a distributed manner  \cite{murray2021partially,Hijazi2012}. In other words, most MINLP algorithms such as ECP, ESH, OA, and GBD are centralized algorithms in the context of distributed computations \cite{Kronqvist2020}.
\end{remark}  
\subsection{A Tour of RH-ADMM}
 Proposed by \cite{olama2019}, the RH-ADMM algorithm is a fully decentralized distributed algorithm to solve general distributed convex optimization problems. Generally speaking, the RH-ADMM algorithm uses the concept of \textit{hypergraphs} to model the CN. A hypergraph is a generalization of an ordinary graph in which an edge can join an arbitrary number of vertices. By following this particular structure, multiple sources, called \textit{Local Fusion Center (LFC)}, can be considered on the CN. This architecture allows decentralized computations over the network, whereby each node can communicate with its corresponding LFC. Besides, the RH-ADMM algorithm benefits from adaptive parameter tuning and a relaxation scheme that improves practical convergence. In synthesis, the main motivations to use RH-ADMM are:
\begin{itemize}
    \item The RH-ADMM algorithm combines centralized and distributed optimization methods to scale the solution of convex problems over large-scale networks. This scaling is achieved by distributing local fusion centers throughout the network.
    
    \item By deriving RH-ADMM from \textcolor{black}{the operator and fixed-point theory \cite{ryu2016primer, olama2019}}, the algorithm is augmented with a relaxation parameter that considerably affects its practical convergence.

    \item RH-ADMM provides an adaptive strategy for tuning the penalty parameter on (primal and dual) residuals that improves convergence and numerical stability.

    \item Finally, the numerical experiments reported by \cite{olama2019} demonstrate the effectiveness of RH-ADMM for solving distributed convex problems.     
    
\end{itemize}
\section{MINLP Reformulation and Solution Architecture}\label{reformlulation}

In this section, MINLP reformulations and a numerical architecture for the SCP problem \eqref{dis-ccp} are presented.
\subsection{Equivalent MINLP Formulation for SCP}
By following the main idea of RH-ADMM, we first introduce multiple LFC nodes to the CN. Consider a hypergraph $\mathcal{H}=(\mathcal{V},\mathcal{E})$ defined as follows: $\mathcal{V}=\{1,2, ...,N\}$ is the set of nodes such that node $i$ decides upon the values of vector variable $\mathbf{x}_i$; $\mathcal{E}=\{\mathcal{E}_k\subset \mathcal{V} : k=1,\ldots,K\}$ is the set of hyperedges, where a hyperedge  $\mathcal{E}_k$ connects all nodes $i\in\mathcal{E}_k$ and $K$ is the number of hyperedges.
   Now we introduce the concept of path in a hypergraph: a path $p(i,j)=\langle \mathcal{E}'_1,..., \mathcal{E}'_k\rangle$ connects nodes $i$ and $j$ if $i\in \mathcal{E}'_1$, $j\in \mathcal{E}'_k$, and $\mathcal{E}'_l\cap \mathcal{E}'_{l+1}\neq \emptyset$ for $l=1,\ldots,k-1$, and $\mathcal{E}'_l\in\mathcal{E}$ for all $l$. 
Put in other words, a path $p(i,j)$ establishes a communication channel between nodes $i$ and $j$, through the hyperedges.
   We make the assumption that the hypergraph $\mathcal{H}$ is connected, meaning that for all $i,j\in\mathcal{V}$ there exists a path $p(i,j)$ connecting $i$ and $j$. Taking into account the hypergraph structure on the CN, the following equivalent formulation for the SCP problem \eqref{dis-ccp} is obtained:
\begin{subequations}\label{admm-ccp}
\begin{align}
\text{P$_2$:}\quad    \min_{\textcolor{black}{\Tilde{\x}, \Tilde{\y}}} &\sum_{i=1}^{N}f_i(\x_i) \label{admm-ccp:eq2a}\\
  s.t.~  &g_h(\x_i) \leq 0, \, \forall i=1,...,N, \, h=1,...,m\\
    &\x_i = \y_j, \, \forall i \in \mathcal{E}_j, \, \mathcal{E}_j \in \mathcal{E} \label{admm-ccp:eq2c} \\
    & \x_i \in \Omega_i, \, \forall i=1,...,N \\
    &\norm{\y_j}_0 \leq \kappa, \, \forall j=1,\ldots,K .
    \label{card}
\end{align}
\end{subequations}
where $\x_i,\y_i\in\mathcal{R}^n$, $\Omega_i = \{\x_i \in \mathcal{R}^n : \mathbf{D}\mathbf{x}_i \leq \mathbf{d}, \mathbf{A}\mathbf{x}_i = \mathbf{b} \}$ is a polytope\footnote{A polytope is a bounded polyhedron \cite[Chapter 2, pg. 31]{boyd2004convex}.}, \textcolor{black}{$\Tilde{\x} = [\x_1^T,...,\x_N^T]^T \in\mathcal{R}^{Nn}$ is the vector of decision variables, and $\Tilde{\y}= [\y_1^T,...,\y_K^T]^T\in\mathcal{R}^{Kn}$  is the vector of auxiliary variables (related to LFCs) which are associated hyperedges.} \textcolor{black}{Notice that, from the assumption that the hypergraph $\mathcal{H}$ is connected, constraint \eqref{admm-ccp:eq2c} is equivalent to $\x_1=\x_2=\cdots=\x_N$.} In particular, here $\x_i$ is a local copy of the global variable $\mathbf{x}$ in $i$-th node and $\mathbf{y}_j$ is the $j$-th LFC variable.
Note that, in contrast to \eqref{D-CCP}, \eqref{dis-ccp} is coupled through the equality constraints \eqref{admm-ccp:eq2c}. To handle such constraints in the primal problem of the OA, we decompose the problem into $N$ separable subproblems by utilizing the RH-ADMM algorithm \cite{olama2019}. Moreover, we impose the sparsity constraint for each LFC to improve the performance of the primal solver to reach a consensus among all nodes in the network. In this case, the primal solver converges in a fewer number of iterations since all LFCs satisfy constraints \eqref{card}.

Problem \eqref{admm-ccp} without the sparsity constraint can be solved in a fully decentralized manner with the recently proposed RH-ADMM algorithm \cite{olama2019}. By introducing suitable binary variables and using the big-M method, Problem \eqref{admm-ccp}  is recast as:
\begin{subequations}\label{admm-minlp}
\begin{align}
\text{P$_3:$}\quad    \min_{\Tilde{\x}, \Tilde{\y}, \Tilde{\zz}} &\sum_{i=1}^{N}f_i(\x_i) \label{admm-minlp:eq0} \\
   s.t.~ &g_h(\x_i) \leq 0, \, \forall i=1,...,N, \,h=1,...,m \label{admm-minlp:eq1 }\\
    &\x_i = \y_j, \, \forall i \in \mathcal{E}_j, \, \mathcal{E}_j \in \mathcal{E} \label{admm-minlp:eq2} \\
    & \x_i \in \Omega_i, \, \forall i=1,...,N \label{admm-minlp:eq3} \\
    &-M_j\zz_j \leq \y_j \leq M_j \zz_j, \forall j=1,...,K \label{admm-minlp:eq4} \\
    &e_j^T\zz_j \leq \kappa, \forall j=1,...,K \label{admm-minlp:eq5}
\end{align}
\end{subequations}
\textcolor{black}{in which $\Tilde{\zz} = [\zz_1^T,...,\zz_K^T]^T\in\mathcal{R}^{Kn}$ is the vector of binary variables, where $\zz_j \in \{0,1\}^{n}, \forall j = \{1,...,K\}$,} and $M_j$ is a sufficiently large constant that can be computed offline, given that the $\Omega_i$ sets are polytopes and therefore bounded. The big-M value $M_j$ can be obtained from data in statistical learning problems \cite{Bertsimas2017} and from the physical bounds in control applications \cite{Aguilera2017}.

Note that since in problem \eqref{admm-minlp} we define the sparsity constraint on each local auxiliary variable $\y_j$, instead of each local $\x_i$ a small number of binary variables is needed as $K \ll N$ in practice. In general, Problem \eqref{admm-minlp} is a distributed MINLP equivalent to Problem \eqref{admm-ccp}. 
This MINLP problem is the final reformulation of the SCP problem for which the DiPOA algorithm is developed.
To solve Problem \eqref{admm-minlp}, we incorporate the RH-ADMM algorithm inside the OA method to act as the numerical kernel, which generates the OA optimality and feasibility cuts. By using the DiPOA algorithm, most of the computational burden is handled in a decentralized fashion, which significantly reduces the CPU time of the solution methodology. Moreover, some heuristics and improvements to control the quality and number of OA cuts are proposed.

\subsection{Solution Architecture}\label{hierarchy}

\begin{figure}
	\centering
	\includegraphics[width=0.5\textwidth]{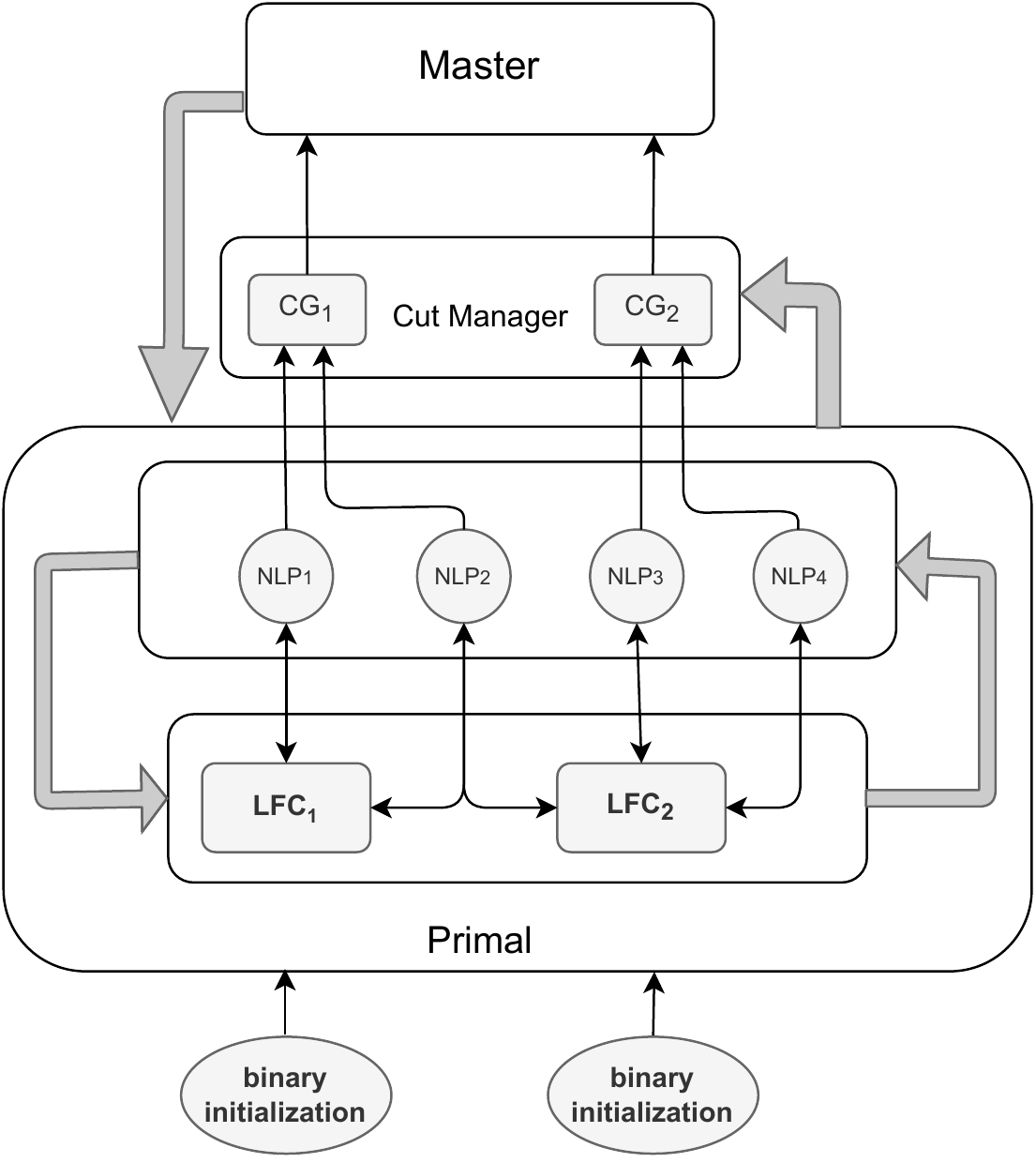}
	\caption{The proposed numerical hierarchy for solving the distributed MINLP Problem \eqref{admm-minlp}.} 
	\label{num-hier}
\end{figure}

Before discussing the development of the DiPOA algorithm, we introduce a multi-level numerical architecture for solving Problem \eqref{admm-minlp}. 
The architecture consists of three main computational levels, namely, \textit{Primal}, \textit{Cut Manager}, and \textit{Master} levels which are depicted in Figure \ref{num-hier}. This figure considers the case where $K = 2$ and $N = 4$. In the following, we explain each level in detail.
\subsubsection*{Primal Level:} 
The primal level deals with the NLP part of Problem \eqref{admm-minlp} and consists of two sub-levels. Each sub-level is also responsible for a certain computational task. To solve the primal NLP part of \eqref{admm-minlp}, we decompose the main computational parts into two steps that distributedly communicate over the CN. In the first sub-level, multiple \textit{local} convex NLP problems are solved in parallel, in a fully decentralized fashion. The local nodes then \textit{synchronously} communicate to the second sub-level which is responsible for another phase of computations. The computations of this level are performed by the LFCs, which consist of the solution of a series of unconstrained NLP problems.
In some practical cases, LFC level computations are related to computing the projection operator, which guarantees the problem feasibility \cite{ma2018fast,olama2019}.
As well as before, the computations in this level are performed in parallel. In many real-world applications, due to the structure of the problem to be solved, the computations at the LFC level are translated into numerical linear algebra operations. 
Therefore, the LFC level is in a tight connection with the LA level and can benefit from parallel LA computations. To guarantee the consensus, at least one node has to be shared between LFCs. This fact is observable in Figure \ref{num-hier} where, for example, NLP$_{2}$ sends its information to both LFC$_1$ and LFC$_2$.
In synthesis, the primal level aims to obtain the solutions to different nonlinear optimization problems\footnote{The problems are \textit{different} in the sense that the problem data are not necessarily the same.} for which the \textcolor{black}{computations are performed in parallel, using a fully decentralized approach.}
\subsubsection*{Cut Manager Level: }
In this level, OA underestimators (\textit{e.g.}, OA cuts) are generated and sent to the master level. In principle, the cut manager needs at least three types of information to generate the OA cuts, which consist of the objective value, the minimizer, and the gradients of each local NLP sub-problem. In particular, the data transfer between LFC and NLP levels continues until the local NLPs reach a consensus and their solutions are optimally feasible. After the consensus, the LFCs send the required information to the CG level where the cuts are generated. This procedure is illustrated with wide gray arrows in  Figure \ref{num-hier}.
\subsubsection*{Master Level: }
This level corresponds to solving the DiPOA master MILP problem based on the cuts generated by the cut manager. The master MILP problem approximates the nonlinear functions of Problem \eqref{admm-minlp}. As the number of cuts increases, this approximation improves until a good piecewise outer approximator is achieved. The binary solution of the master level is then sent to the LFC level which, together with the NLP level, solves another set of NLP problems. 

Up to this point, we discussed the numerical levels of the proposed solution hierarchy for Problem \eqref{admm-minlp}. In the following, we explain the development of DiPOA in more detail.

\section{\textbf{Distributed Primal Outer Approximation (DiPOA) Algorithm}} \label{dipoa}

This section presents the development of the proposed DiPOA algorithm for solving the distributed MINLP \eqref{admm-minlp} which is equivalent to the SCP problem. Moreover, to improve the algorithm efficiency some techniques and heuristics are introduced.

\subsection{General Structure of DiPOA}
We start by reformulating the problem \eqref{admm-minlp} as the epigraph minimization problem which is also called \textit{lifted formulation} \cite{kronqvist2018reformulations}. In the context of DiPOA, we use lifted formulation for each individual objective function. \textcolor{black}{As the main advantage, the lifted formulation will result in tighter outer approximations when linearizing local objective functions \cite{kronqvist2018reformulations}.} By using the epigraph form, the objective functions become linear and the problem nonlinearities only appear in the constraints. Accordingly, we define $N$ upper bounds $\alpha_i\in \mathcal{R}, \, i = \{1,...,N\}$, one for each of local function, $f_i$, and rewrite  problem \eqref{admm-minlp} as follows:
\begin{subequations}\label{admm-minlp-epi}
\begin{align}
\text{P$_4$:}\quad    \min_{\boldsymbol{\alpha}, \Tilde{\x}, \Tilde{\y}, \Tilde{\zz}} &\sum_{i=1}^{N}\alpha_i\\
  s.t.~  & f_i(\x_i) \leq \alpha_i, \, i = 1,..., N\\
    &g_h(\x) \leq 0, \,\, \forall h=1,...,m\\
    &\x_i = \y_j, \, \forall i \in \mathcal{E}_j, \, \mathcal{E}_j \in \mathcal{E}, \, \forall j = 1,...,K\\
    & \x_i \in \Omega_i\\
    &-M_j\zz_j \leq \y_j \leq M_j \zz_j\\
    &e_j^T\zz_j \leq \kappa
\end{align}
\end{subequations}
where $\boldsymbol{\alpha} = [\alpha_1,...,\alpha_N]^T$. Note that since the local $f_i$'s are convex functions, their epigraph is also convex. Hence, the convexity of problem \eqref{admm-minlp-epi} is preserved and problems \eqref{admm-minlp} and \eqref{admm-minlp-epi} are equivalent. 

According to the first-order convexity property of convex functions, each local function $f_i$ can be globally underestimated by a linear supporting hyper-plane around the point $\bar{\x}_i$. This approximation can be written as:
\begin{equation}\label{lin-app}
    f_i(\x_i) \geq f_i(\bar{\x}_i) + \nabla f_i(\bar{\x}_i)^T (\x_i - \bar{\x}_i), \,\, \forall i = 1,..., N
\end{equation}
Moreover, as $f_i(\x_i) \leq \alpha_i$, inequality \eqref{lin-app} can be written as,
\begin{equation}\label{lin-app-epi}
    \alpha_i\geq f_i(\bar{\x}_i) + \nabla f_i(\bar{\x}_i)^T (\x_i - \bar{\x}_i), \,\, \forall i = 1,..., N
\end{equation}
  By a similar argument and using the fact that  $g(\x)_i \leq 0$, we have
\begin{equation}\label{lin-app-feasb}
    g_h(\Bar{\x}_i) + \nabla g_h(\Bar{\x}_i)^T (\x_i - \Bar{\x}_i) \leq 0, \,\, \forall i = 1,..., N
\end{equation}
  
The inequalities provided by equations \eqref{lin-app-epi} and \eqref{lin-app-feasb} are referred to as \textcolor{black}{OA cuts \cite{duran1986outer}}. Considering all feasible points, equivalently, problem \eqref{admm-minlp-epi} can be stated as the following distributed Mixed-Integer Linear Programming (D-MILP) problem:
\begin{subequations} \label{master-milp}
\begin{align}
\text{P$_5$: (D-MILP})\quad    \min_{\boldsymbol{\alpha}, \Tilde{\x}, \Tilde{\y}, \Tilde{\zz}} &\sum_{i=1}^{N}\alpha_i\\
  s.t.~  & \alpha_i\geq f_i(\bar{\x}_i) + \nabla f_i(\bar{\x}_i)^T (\x_i - \bar{\x}_i), \,\, \forall \bar{\x}_i \in \mathcal{X}_i \label{eq:P7:optimality-cut} \\
    &0 \geq  \nabla g_h(\bar{\x}_i)^T (\x_i - \bar{\x}_i), \,\, \forall \bar{\x}_i \in \partial\mathcal{X}_i, h\in\mathcal{A}(\bar{\x}_i) \label{eq:P7:feasibility-cut} \\
    &\x_i = \y_j, \, \forall i \in \mathcal{E}_j, \, \mathcal{E}_j \in \mathcal{E}, \, \forall j = 1,...,K\\
    & \x_i \in \Omega_i\\
    &-M_j\zz_j \leq \y_j \leq M_j \zz_j\\
    &e_j^T\zz_j \leq \kappa
\end{align}
\end{subequations}
where:
\begin{itemize}

\item $\mathcal{X}_i = \{\x_i: g(\x_i) \leq 0, \, \x_i \in \Omega_i\}$ is the set of all feasible points. Notice that $\mathcal{X}_i = \mathcal{X}$ for all $i$, but we keep the index $i$ to make it clear that the problem refers to node $i$.

\item $\partial\mathcal{X}_i$ is the set of feasible points that lie on the boundary of the feasible set.

\item $g(\mathbf{x}) = [g_1(\mathbf{x}), ...,g_m(\mathbf{x})]^T$ is a vector function with the nonlinear constraints.

\item $\mathcal{A}(\bar{\x}_i) = \{h : g_h(\bar{\x}_i)=0\}$ is the index set of active nonlinear constraint functions at point $\bar{\x}_i$.
\end{itemize}
The outer approximations \eqref{eq:P7:optimality-cut} are obtained from the linearization of the objective function for all points $\bar{\x}_i$ that are feasible with respect to $\Omega_i$ and the nonlinear constraints (the points in $\mathcal{X}_i$). On the other hand, the outer approximations of the constraints \eqref{eq:P7:feasibility-cut} are obtained by linearization of only the \textit{active} nonlinear constraint functions in $g(\cdot)$, at the points on the boundary of $\mathcal{X}_i$, namely the set $\partial \mathcal{X}_i$.

\textcolor{black}{According to \cite{duran1986outer},  the  MILP problem \eqref{master-milp} involves an infinite number of constraints to be equivalent to the MINLP problem \eqref{admm-minlp-epi}. Instead, we follow the strategy of the outer approximation algorithm that iteratively solves a relaxation of problem \eqref{master-milp} and adds OA cuts on demand. Following this strategy, the following relaxation of \eqref{master-milp} is used as:}
\begin{subequations}\label{relaxed-master-milp}
\begin{align}
\text{P$_6$: }\quad    \min_{\alpha, \Tilde{\x}, \Tilde{\y}, \Tilde{\zz}} &\sum_{i=1}^{N}\alpha_i\\
  s.t.~  & \alpha_i\geq f_i(\x_i^k) + \nabla f_i(\x_i)^T (\x_i - \x_i^k), \,\, \forall \x^k_i \in \mathcal{X}_i^k  \\
    &0 \geq  \nabla g_h(\x_i^k)^T (\x_i - \x_i^k), \,\, \forall \x_i^k\in \partial\mathcal{X}_i^k, h\in\mathcal{A}^k(\x_i^k) \\
    &\x_i = \y_j, \, \forall i \in \mathcal{E}_j, \, \mathcal{E}_j \in \mathcal{E}, \, \forall j = 1,...,K\\
    & \x_i \in \Omega_i\\
    &-M_j\zz_j \leq \y_j \leq M_j \zz_j\\
    &e_j^T\zz_j \leq \kappa
\end{align}
\end{subequations}
where 
$k$ denotes the iteration counter and $\mathcal{X}_i^k$ and $\mathcal{A}_i^k$ are finite sets consisting of local feasible points defined as 
\begin{align*}   
\mathcal{X}_i^k &= \{\x_i^\ell: g_h(\x_i^\ell) \leq 0, \, \x_i^\ell \in \Omega_i, \, \forall \ell \in \{1,\ldots,k\}\}, \\
\mathcal{A}(\x_i^k) &= \{j : g_j(\x_i^k)=0\}, 
\end{align*}
Particularly, $\mathcal{X}_i^k$ consists of the feasible points up to the current iteration $k$ while $\mathcal{A}_i^k$ includes the indices $j$ for which the nonlinear constraints are active. We refer to problem \eqref{relaxed-master-milp} as the DiPOA \textit{master's problem}. Since \eqref{relaxed-master-milp} is a relaxation, its optimal value provides a valid lower bound on the optimal objective value. 

Typically, there are various ways to obtain local feasible points, $x_i^k$'s, and generate OA and feasibility cuts. For instance, the cuts can be generated based on the fractional solution of the relaxed problem \eqref{relaxed-master-milp}, where the integrality constraint on binary variables is dropped. Although this method is extremely fast to generate the cuts, the points about which the cuts are taken may be far from feasible, so the resulting linearizations may form a poor approximation of the nonlinear feasible region. An efficient way to generate feasible solutions for each node of the CN can be achieved by fixing the local binary decision variables and solving the resulting optimization problem, which is obtained as follows:
\begin{subequations}\label{dist-prim}
\begin{align}
\text{P$_7$}(\zz_1^k,...,\zz_K^k)\quad   \min_{\Tilde{\x}, \Tilde{\y}} &\sum_{i=1}^{N}f_i(\x_i) \\
 s.t.~  & g_h(\x_i) \leq 0, \,\, \forall h=1,...,m \label{nlc}\\
        &\x_i = \y_j, \, \forall i \in \mathcal{E}_j, \,\, \mathcal{E}_j \in \mathcal{E}\\
        & \x_i \in \Omega_i, \,\, \forall i = 1,..., N\\
        & -M_j \zz_j^k \leq \y_j \leq M_j \zz_j^k,\,\, \forall j = 1,..., K
\end{align}
\end{subequations}
where $\zz_j^k, \forall j = 1,...,K$ are fixed binary variables. We refer to problem \eqref{dist-prim} as the DiPOA \textit{primal problem}. The solution of problem \eqref{dist-prim} has the advantage of generating linearizations about points that are closer to the feasible region.
The optimal value of problem \eqref{dist-prim} yields valid upper bounds on the local $f_i$ functions and their solutions provide the necessary information \textcolor{black}{to generate DiPOA feasibility cuts and outer approximations cuts.}

Problem \eqref{dist-prim} has a suitable form for which the RH-ADMM algorithm can be applied. In the general case, to compute the LFC variables $\y_j$, a bounded constrained optimization problem has to be solved. In this case, each LFC can compute $\y_j$ satisfying  the bounded constraints defined by $M_j$ \cite{olama2019}.

In summary, DiPOA is usually initiated by given feasible binary variables which are utilized to form the primal problem \eqref{dist-prim}. The optimal value of \eqref{dist-prim} provides an upper bound on the optimal value of each local function. Moreover, the optimal solution of problem \eqref{dist-prim} is used to construct the outer approximation supporting hyper-planes. Then, by solving problem \eqref{relaxed-master-milp} a lower bound on the optimal value of each local function $f_i$ is generated. The feasible binary solutions $\zz_i^k$'s of problem \eqref{relaxed-master-milp} are then considered as the fixed integer solutions for the primal problem \eqref{dist-prim} in the next iteration.  Since at each iteration $k$, a limited number of cuts (usually one cut per iteration) is generated and added to the DiPOA master problem,  the outer approximation of the local functions becomes tighter and the generated lower bound increases as the DiPOA progresses. This procedure is repeated until the gap between the upper and lower is within a given tolerance.

\begin{remark}
 In some specific cases (e.g., distributed), solving the NLP sub-problems can be expensive and stands as a bottleneck for the MINLP algorithms (\textit{e.g.}, OA). Taking this issue into account, it can be noticed that the primal problem \eqref{dist-prim} is a consensus distributed non-linear convex optimization problem. To improve computational efficiency, this problem should be solved distributedly. Accordingly, the state-of-the-art NLP solvers cannot be used directly to solve problem \eqref{dist-prim} since most of the standard NLP solvers are centralized. Hence, it is critical to utilize efficient distributed algorithms to deal with the NLP sub-problems, while exploiting the structure of the problem.
\end{remark}

Before presenting the full algorithmic framework and implementation details, in the next section of this paper, we discuss some effective heuristics and improvements that are applied to the SCP problem \eqref{dis-ccp}. In the following, we develop specialized heuristics and improvements to the problem \eqref{dis-ccp} and its equivalent MINLP reformulation to improve the efficiency of the DiPOA algorithm from both convergence and computation aspects. These improvements consist of a \textit{Specialized Feasibility Pump (SFP)} and \textit{Second order Cuts (SoCut)} strategy.
\subsection{Specialized Feasibility Pump}
In general MINLP algorithms, the Feasibility Pump (FP) methods are algorithms that quickly find initial feasible solutions. 
  \textcolor{black}{The FP is a primal heuristic for mixed-integer programming based on the idea of producing two sequences of points that hopefully converge to a feasible solution to a given optimization problem.} 
One sequence consists of points that are feasible for a continuous relaxation, however, possibly integer infeasible. The other sequence consists of integral points that might violate some of the imposed constraints. The next point of one sequence is generated by minimizing the distance to the latest point of the other sequence.

In this section, we propose a Specialized FP (SFP) method for problem \eqref{dis-ccp} which provides a good feasible solution along with a high-quality upper bound. Our approach applies to problem \eqref{dis-ccp} with any convex functions subject to convex and sparsity constraints. The SFP method consists of two main steps:
\begin{enumerate}

    \item \textbf{Relaxation step}: In this step, the RH-ADMM algorithm solves a relaxed version of the (SCP) problem with a convex approximation of the $\ell_0$ norm. Considering the sparsity constraint \eqref{card} and its reformulation \eqref{admm-minlp:eq5} in terms of the LFC variable $y_j$, this restriction can be conveniently expressed by the following set
\begin{equation}
    \mathcal{D}_j = \left\lbrace \y_j\in \mathcal{R}^n\,|\, \norm{\y_j}_0 \leq \kappa \right\rbrace,
\end{equation}
because the bound $M_j$ is valid for $\y_j$, given that the $\x_i$ variables belong to polytopes, this bound can be introduced to $\mathcal{D}_j$ in the definition of the following set:
\begin{equation}
    \mathcal{C}_j = \left\lbrace \y_j\in \mathcal{R}^n\,|\, \norm{\y_j}_0 \leq \kappa \, ,\norm{\y_j}_\infty \leq M_j \right\rbrace . 
\end{equation}
\cite{Bertsimas2015} showed that the convex hull $\mathcal{C}_j$ lies inside a simple set:
\begin{equation}
    \text{Conv}\left( \mathcal{C}_j \right) \subseteq \left\lbrace \y_j\in \mathcal{R}^n\,|\,\norm{\y_j}_1 \leq \kappa M_j \right\rbrace . \label{eq:relaxed-Conv:C_j}
\end{equation}
Hence, the problem which is solved in the relaxation step of SPF is as follows,
\eqref{dis-ccp} 
\begin{subequations}\label{dist-prim-sfp}
\begin{align}
\qquad \min_{\Tilde{\x}, \Tilde{\y}} &\sum_{i=1}^{N}f_i(\x_i) \\
 s.t.~   & g_h(\x_i) \leq 0, \,\, \forall h=1,...,m\\
    &\x_i = \y_j, \, \forall i \in \mathcal{E}_j, \,\, \mathcal{E}_j \in \mathcal{E}, \,\, \forall j = 1,..., K\\
    & \x_i \in \Omega_i, \,\, \forall i = 1,..., N\\
    & \y_j \in \left\lbrace \y_j\in \mathcal{R}^n\,|\,\norm{\y_j}_1 \leq \kappa M_j \right\rbrace,\,\, \forall j = 1,..., K
\end{align}
\end{subequations}

Problem \eqref{dist-prim-sfp} is a distributed convex NLP for which RH-ADMM can be applied. Since the sparsity constraint is relaxed according to \eqref{eq:relaxed-Conv:C_j}, the optimal values of \eqref{dist-prim-sfp} provide a valid lower bound for SCP.
Therefore, the relaxation step of the SFP method consists of solving problem \eqref{dist-prim-sfp} using the RH-ADMM algorithm. Accordingly, at each iteration of the RH-ADMM algorithm, the local solutions, $\x_i$, are updated by each node of the CN and sent to the LFC level where the relaxed 1-norm constraint \eqref{eq:relaxed-Conv:C_j} is computed. This procedure is repeated until RH-ADMM converges.
    
    \item \textbf{Projection step}: This step consists of projecting the local solutions provided by the first step of the SFP to the feasible region of the original sparsity constrained problem.  The feasible solution provided by this step defines a good initial upper bound on the optimal objective values. Moreover, the obtained local feasible solutions are used to approximate the local objective functions by generating cuts. These cuts provide a good linearization point for the DiPOA algorithm and help the algorithm to converge much more quickly. A similar approach is used in \cite{Bernal2020} to improve the performance of the solver DICOPT in convex MINLP problems.
   In the following, we provide details regarding the SFP steps.  Before discussing the projection step of the SFP method, we define the projection on a set.

\begin{definition}
The projection of a point $\y \in \mathcal{R}^n$ onto a closed set $\mathcal{A}\subseteq\mathcal{R}^n$ is denoted by $\y^*$ and defined as:
\begin{equation}\label{projection-step}
    \y^* = \Pi_{\mathcal{A}}(\y) = \arg\min_{\textcolor{black}{\y}\in\mathcal{A}} \norm{\y_0 - \y} 
\end{equation}
\end{definition}

In some special cases, the projection of a point on a set is unique. For instance, if $\mathcal{A}$ is a closed and convex set and the norm is strongly convex, then the projection is unique. However, for non-convex sets (\textit{e.g.}, $\mathcal{C}_j$) the projection is not unique and can be difficult to achieve. However, for some specific cases, it is not too difficult to compute the projection on a non-convex set. In the case of the set $\mathcal{C}_j$ in which the zero-norm is incorporated, we extend the definition of $\Pi_{\mathcal{C}_j}(\y)$  from  \cite{boyd2011distributed} with a tie-breaking rule:
\begin{definition}
The projection $\Pi_{\mathcal{C}_j}(\y)$ of a point $\y\in \mathcal{R}^n$ onto the set $\mathcal{C}_j$ keeps the $\kappa$ largest (in absolute value) elements and zeros out the remaining, breaking ties in the lexicographic order. 
\end{definition}
\textbf{Example 1.} Consider $\y \in \mathcal{R}^5$ to be defined as $\y = \begin{pmatrix} 5.7& 1.4& -3.2& -2.3 & 2.3   \end{pmatrix}^T$ and $\kappa = 3$. The projection is obtained as follows:
\begin{equation*}
   \y^* = \Pi_{\mathcal{C}}(\y) = \begin{pmatrix} 5.7& 0& -3.2& -2.3& 0  \end{pmatrix}^T
\end{equation*}
Based on $\y^*$ we can generate a binary vector $\zz \in \{0, 1\}^5$, such as
$\mathbf{z} = \begin{pmatrix} 1& 0& 1& 1& 0  \end{pmatrix}^T$ in the example, that can be used to initialize the DiPOA algorithm. The tie between $-2.3$ and $2.3$ is resolved by selecting the element that comes first in the lexicographic order.

Therefore, the projection vector, $\y^*$, would be a feasible solution for SCP, which provides an upper bound on the original problem and can be used to start DiPOA. If $\y^*$ is not feasible, then we can attempt to recover feasibility by reoptimizing the primal's problem with RH-ADMM, but allowing to be nonzero only the entries of $\y_j$ that are projected onto $\mathcal{C}_j$, thereby ensuring $\|\y_j\|_0\leq  \kappa$.
\end{enumerate}
\subsection{Lower Bound Improvement}
As mentioned, at each iteration of decomposition-based MINLP algorithms (\textit{e.g.}, OA, GBD, and ECP), a linear cut is generated and added to the master's problem to approximate the nonlinear objective and constraint functions. By enforcing the cuts at each iteration and solving the master MILP problem, a tighter approximation is achieved at each iteration. However, in some cases (\textit{e.g.}, when highly nonlinear functions are incorporated into the optimization problem), the linear approximation can be poor. In this case, the lower bound provided by solving the MILP master problem will be improved slowly and a large number of outer approximation iterations, $k$, is needed. To improve the quality of the approximation cuts,  several methods have been introduced in the literature. Among all methods, two well-known methods are:
\begin{enumerate}

    \item \textbf{Multi-Generation Cut methods.} Originally, a single cut per iteration is generated for the MILP problem of the OA algorithm \cite{su2015computational}. To improve the approximation of the nonlinear functions, several cuts can be added in each iteration. One possibility to generate multiple cuts for the master's problem is to use the solution pool of the MILP solver (e.g., CPLEX and GUROBI), as implemented by the solver SHOT \cite{Lundell2019}.
     By using these feasible points, multiple instances of the primal problem can be solved in parallel to yield multiple cuts, often improving overall performance.  

    \item \textbf{Second order Cut (SoCut) methods.} In these methods, the second-order information of the objective and constraint functions are used to formulate the master's problem. There are two well-known methods to design the master's problem using the SoCuts. 
    In the first method,  the master problem is built based on the Hessian information of the Lagrange function.
    In this case, the master's problem consists of a quadratic objective function and linear constraints which results in a Mixed-Integer Quadratic Programming (MIQP) problem \cite{Kronqvist2020}. In the second method, the second-order information of the nonlinear functions is used to generate quadratic underestimators. Quadratic outer approximations proved to be a stricter and tighter underestimation for convex nonlinear functions than classical supporting hyper-planes. Although  SoCuts increase the complexity of the master's problem, they can provide significant improvements in generating valid lower bounds.
\end{enumerate}
Regarding the second class of SoCut-based methods,  \cite{su2018improved} showed that using scaled Taylor series truncated at the second-order and with a scaling parameter $\beta \in \{0,1\}$, a valid underestimator for convex nonlinear functions is generated. The scaling parameter $\beta$ is introduced since the unscaled second-order Taylor approximation is not tight and global in general and, hence, by choosing a sufficiently small $\beta$, a global underestimator is achieved.  However, since selecting the optimal $\beta$ is not easy, an online enumeration algorithm based on vertex enumeration was proposed. 
\subsubsection{Second Order Approximation as a Global Underestimator:}

In many practical applications of the SCP problem (\textit{e.g.}, learning and control), the nonlinear functions are not only convex but also strongly convex functions. By assuming the strong convexity property on the nonlinear functions, we specialize the method of \cite{su2018improved} in such a way that a distributed global underestimator is achieved and no online tuning is needed. Utilizing this framework, the tightest global quadratic approximation of the local nonlinear functions is derived. 

We also introduce an event-triggered framework to generate and add the  SoCuts, which is herein referred to as the Event-Triggered Second-Order Cut (ET-SoCut) method. 
Before discussing the ET-SoCut procedure, we provide the definition of $m$-strong convexity of nonlinear functions.
\begin{definition}
A twice differentiable function $f: \mathcal{R}^n\xrightarrow{}\mathcal{R}$ is $m$-strongly convex if there exists an $m> 0$ such that 
\begin{equation}\label{sconv}
    \nabla^2f(\x) \succeq m \qquad \forall \x \in \text{dom}\,f
\end{equation}
where $\nabla^2f(\x) \in \mathcal{R}^{n\times n}$ is the Hessian matrix at $\x$.
\end{definition}

\begin{figure}[htb]
    \begin{subfigure}{0.5\textwidth}
		\includegraphics[width=0.85\textwidth]{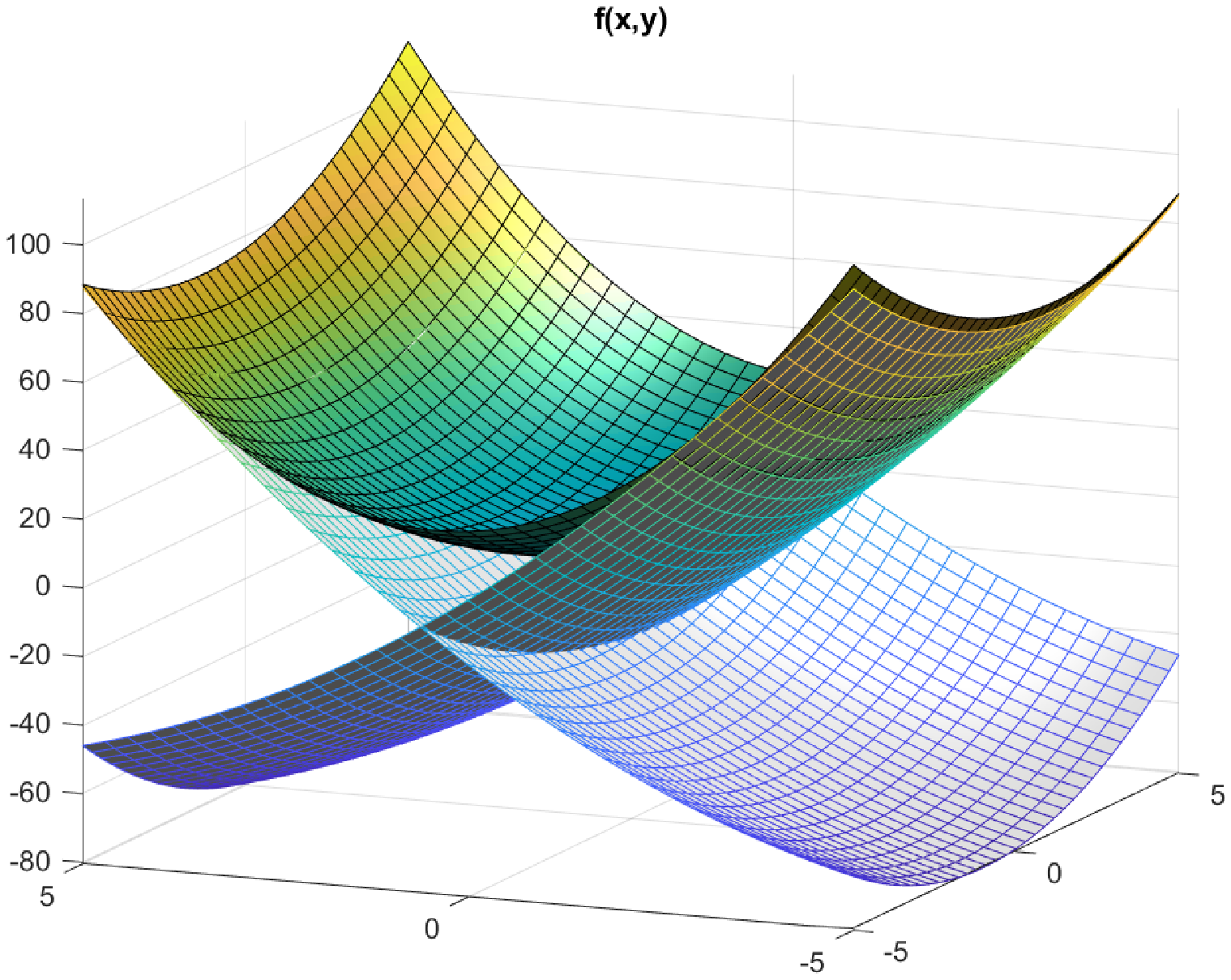}
		\caption{}
		\label{app-a}
    \end{subfigure}
         \hfill
        \begin{subfigure}{0.5\textwidth}
	\includegraphics[width=0.85\textwidth]{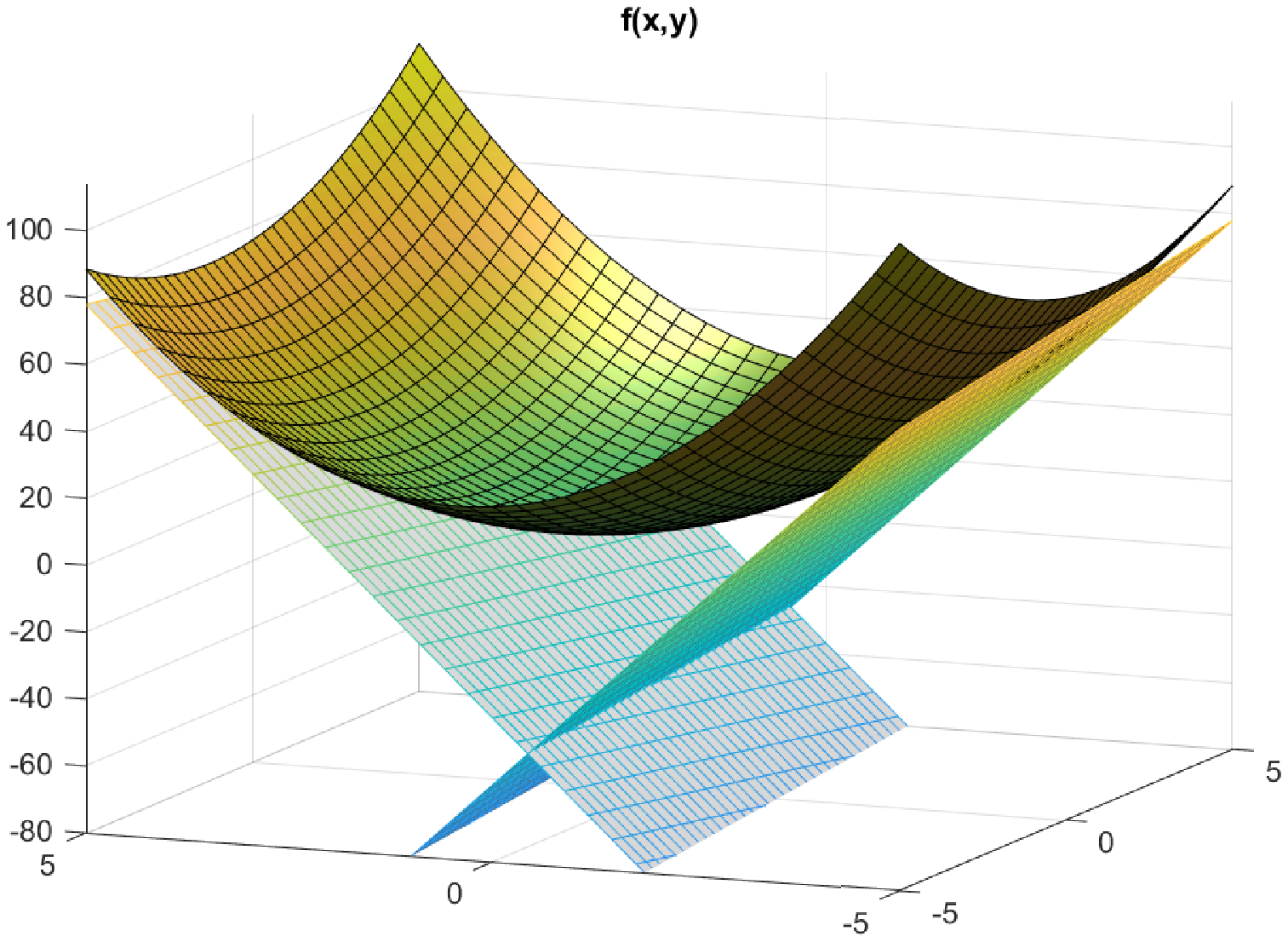}
		    \caption{}
		    \label{app-b}
    \end{subfigure}
    \caption{Comparison between linear and quadratic underestimators. a) quadratic underestimators. b) linear underestimators} 
    \label{cut-pic}
	\end{figure}

By assuming strong convexity of the local functions $f_i$ associated to the CNs, a global quadratic underestimator can be achieved for each $f_i(\x_i)$ about a point $\x_i^k$:
\begin{equation}\label{fineq}
     f_i(\x_i) \geq f_i(\x_i^k) + \nabla f_i(\x_i^k)^T (\x_i - \x_i^k) + \frac{m_i}{2}\norm{(\x_i - \x_i^k)}_2^2, \quad \forall i \in \{1,...,N\}.
\end{equation}
Figure \ref{cut-pic} shows an example that illustrates the advantage of using SoCuts over linear OA cuts. The inequality \eqref{fineq} has two important properties:
\begin{enumerate}
    \item The local Hessian matrices do not need to be computed and stored. 
    \item The computation of $m_i$ is straightforward in some cases. 
\end{enumerate} 

For general strongly convex functions, it is more challenging to compute $m$ efficiently. However, in practical problems of learning and control, the computation of $m$ can be straightforward. For example, in many instances of sparse Model Predictive Control (s-MPC) problems \textcolor{black}{(which is a subclass of the SCP problem)} the objective function is a convex quadratic function. In s-MPC, the goal is to control a process employing a reduced number of inputs, which can benefit the operation of the control system \cite{Aguilera2017}.  For convex quadratic functions, $m$ is the smallest Eigenvalue of the Hessian matrix. In machine learning problems the objective function usually consists of a convex function and a strongly convex \textit{regularization} term. In this case, $m$ can be computed from the regularization term. More details about implementing SoCuts are discussed later.

In summary, based on \eqref{fineq}, we generate quadratic cuts for the nonlinear convex functions of the SCP problem \eqref{dis-ccp} to accelerate the convergence of the DiPOA algorithm.
\subsubsection{Global SoCuts for the DiPOA Algorithm: }
By considering the global SoCuts, we can form the following cuts for the master problem of the DiPOA algorithm
\begin{equation}\label{soc}
     \alpha_i \geq f_i(\x_i^k) + \nabla f_i(\x_i^k)^T (\x_i - \x_i^k) + \frac{m_i}{2}\norm{(\x_i - \x_i^k)}_2^2, \quad \forall i \in \{1,...,N\}.
\end{equation}
Hence, the relaxed master's problem becomes,
\begin{subequations}\label{relaxed-master-miqcp}
\begin{align}
\min_{\alpha, \Tilde{\x}, \Tilde{\y}, \Tilde{\zz}} &\sum_{i=1}^{N}\alpha_i\\
  s.t.~  & \alpha_i \geq f_i(\x_i^k) + \nabla f_i(\x_i^k)^T (\x_i - \x_i^k) + \nonumber\\
    &\qquad\frac{m_i}{2}\norm{(\x_i - \x_i^k)}_2^2, \,\forall \x_i^k \in \mathcal{X}_i^k\\
    &0 \geq \nabla g_h(\x_i^k)^T (\x_i - \x_i^k),  \,\, \forall \x_i^k\in \partial\mathcal{X}_i^k, h\in\mathcal{A}(\x_i^k) \\
    &\x_i = \y_j, \, \forall i \in \mathcal{E}_j, \, \mathcal{E}_j \in \mathcal{E}, \, \forall j = 1,...,K.\\
    & \x_i \in \Omega_i\\
    &-M_j\zz_j \leq \y_j \leq M_j \zz_j\\
    &e_j^T\zz_j \leq \kappa
\end{align}
\end{subequations}
Problem \eqref{relaxed-master-miqcp} is a distributed Mixed-Integer Quadratically Constrained Programming (MIQCP) problem. Compared to \eqref{relaxed-master-milp}, since problem \eqref{relaxed-master-miqcp} approximates the nonlinear functions more accurately, a high quality lower bound for the SCP problem is achieved which can help the DiPOA algorithm converge faster. 

\color{black}
The convergence theorem of OA \cite{duran1986outer, Fletcher1996} holds for the DiPOA algorithm because the   generated SoCuts do not cut off the feasible region of the original SCP problem \eqref{admm-ccp}. This property of DiPOA leads to the following remark.
\begin{remark} For a sufficiently large number of outer approximations, $k \geq k^*$, the optimal solution of problem \eqref{admm-ccp} and the master's problem \eqref{relaxed-master-miqcp} are identical.
\end{remark}
\color{black}


 Typically, MIQCP problems are more difficult to solve and require longer CPU time than MILPs. 
However, recently, the capabilities of MIP solvers such as CPLEX and GUROBI have been extended to solve MIQCP problems. Notably, the nature of the Hessian matrix is a crucial factor in solving convex MIQCPs.
 Compared to convex MIQCPs with dense Hessian matrices, general MIP solvers should find it easier to solve problems with constraints in form \eqref{soc} that approximate the Hessian with a quadratic diagonal matrix, and whose elements are the minimum eigenvalue of the Hessian. This is an advantage of the proposed SoCuts  \eqref{soc}, which can equip the master problem \eqref{relaxed-master-miqcp} to be solved more easily when compared to the method of scaled quadratic cuts.

\subsubsection{Event-Triggered Second Order Cut (ET-SoCut)}

 Although MIQCP problems can be solved efficiently with today's solvers, keeping the master's problem simple and easy to solve is important for the efficiency of the DiPOA algorithm.
 To keep the master's problem as simple as possible, and also benefit from better approximations (SoCuts), we propose an Event-Triggered Second-Order Cut (ET-SoCut) generation strategy in which the  SoCuts are generated whenever necessary. To do so, we define an event that is based on the relative optimality gap of the DiPOA algorithm. The relative optimality gap at each DiPOA iteration $k$ is defined as,
 \begin{equation}\label{rgap}
     r^k= \frac{Ub^k - Lb^k}{\max\{Ub^k,0.001\}} \times 100
 \end{equation}
 where $Ub$ and $Lb$ are the upper and lower bound, respectively. 
 Based on the relative gap \eqref{rgap}, at each iteration $k$, we define
 \begin{equation}\label{event}
     e^k = \frac{ r^{k - 1} - r^k   }{r^{k - 1} } 
 \end{equation}
 as an event for the master's problem. At each iteration $k$, $e^{k}$ measures the difference between two consecutive relative gaps. This event is triggered, whenever, 
 \begin{equation}
     e^k \leq \epsilon
 \end{equation}
 where $\epsilon > 0$ is a small number. Therefore, using the  ET-SoCut strategy, at each iteration of the DiPOA algorithm, the relative increment error is computed and if it is smaller than a threshold, a  SoCut is generated for the master's problem. In other words, to avoid generating the  SoCuts at each iteration, we just add  SoCuts whenever the relative gap starts to flatten out. Otherwise, the linear approximations are generated as usual. By using  ET-SoCut, a small number of  SoCuts are generated to help the MIP solver for solving the problem efficiently. More details about the ET-SoCut strategy are provided in the implementation section. 
\subsection{Practical Infeasibility Detection}
This section introduces a practical algorithm to detect the nonlinear infeasibility of problem \eqref{dist-prim} before attempting to solve the problem. Although problem \eqref{dist-prim} provides reliable feasible points at which linear and quadratic approximations are generated, however, certain combinations of fixed binary variables, $\mathbf{z}_j^k$, can lead to the infeasibility of the nonlinear constraint \eqref{nlc}. Therefore, it becomes essential to detect the infeasibility beforehand by recognizing an infeasible point. In such a situation, one solution is to generate a constraint, the so-called \textit{feasibility cut}, that cuts off the infeasible point from the feasible region \cite{Fletcher1996}. The feasibility cuts are usually obtained by solving $\ell_1$ or $\ell_\infty$ norm minimization of infeasible nonlinear constraints. Therefore, it is essential to distinguish the infeasible and feasible nonlinear constraints, which is a task left to the NLP solver by the original OA algorithm. In DiPOA, however, the infeasibility detection feature is not inherently supported.
To detect infeasibility of the local nonlinear constraints, we form the following optimization problem.
\begin{subequations}\label{dist-prim-inf}
\begin{align}
\min_{\Tilde{\x}, \Tilde{\y}} &\sum_{i=1}^{N}\sum_{h=1}^{m}g_h(\x_i) \\
  s.t.~ 
    &\x_i = \y_j, \, \forall i \in \mathcal{E}_j, \,\, \mathcal{E}_j \in \mathcal{E}\\
    & \x_i \in \Omega_i, \,\, \forall i = 1,..., N\\
    & -M_j \zz_j^k \leq \y_j \leq M_j \zz_j^k,\,\, \forall j = 1,..., K \label{bounded}
\end{align}
\end{subequations}
This problem is a consensus convex optimization problem with linear constraints which can be solved using DiPOA. Moreover, since $\Omega_i$ is a polytope and the bounded constraint \eqref{bounded} cannot be empty, problem  \eqref{dist-prim-inf} is always feasible. If the optimal solution of \eqref{dist-prim-inf} yields a positive objective function, then this solution serves a certificate that the corresponding distributed NLP problem \eqref{dist-prim} is infeasible. In this case, the next step is to detect the indices of the infeasible nonlinear constraints. Considering $\Bar{\textbf{x}}_i, i = 1,..., N$, as the optimal solution of \eqref{dist-prim-inf}, we can easily detect the nonlinear constraints by computing $g_h(\Bar{\textbf{x}}_i), i = 1,...,N,\, h=1,...,m$. Finally, we generate and introduce \eqref{lin-app-feasb} around $\Bar{\textbf{x}}_i$ and solve the master's problem \eqref{relaxed-master-miqcp}. In the centralized case ($N = 1)$,  \cite{Fletcher1996} showed that if \eqref{dist-prim} is infeasible, for the given $\mathbf{z}_j^k$, and $\Bar{\textbf{x}}_i$ is an optimal solution of \eqref{dist-prim-inf}, then inequality \eqref{lin-app-feasb} cuts off $\mathbf{z}_j^k$.

\section{Implementation}\label{implementation}
  This section presents the details of the implementation of the DiPOA algorithm.
Before formalizing the full algorithm steps, a summary of parallel programming is introduced. 

\subsection{Parallel Programming and Message Passing Interface (MPI)}
Taking into account the computer hardware structure available today, there exist three levels of parallelism that can be utilized. 
The lowest level of the hierarchy of parallelism regards \textit{vectorization} which supports the \textit{single instruction multiple data} programming style. 
In this case, a limited number of instructions can be performed per CPU cycle in parallel. Most of the compilers of the high-performance programming languages such as {\tt C/C++}  and {\tt FORTRAN}  support vectorization inherently. 
The next level of parallelism consists of \textit{multi-threading}, which allows the programmer to perform \textit{shared memory} parallel computations using different threads on a single machine. 
 Finally, on the highest level of the parallelism hierarchy, the Message Passing Interface (MPI) is implemented \cite{gropp1999using} to allow \textit{distributed computations} on a \textit{distributed memory} machine, such as high-performance computer clusters.
 As the main benefit, MPI supports \textit{data-passing} between various computing nodes, where the local computations are performed. Considering that each computing node in MPI can perform multi-threaded computations, it is possible to use MPI and multi-threaded programming at the same time. This style of programming related to distributed computation is called the \textit{hybrid programming paradigm}. In this paper, we implement the DiPOA algorithm according to the hybrid programming style to take advantage of both parallel computation styles. 

The parallel implementation is one of the most important features of the RH-ADMM algorithm, which is used to solve the D-NLP sub-problems of the DiPOA algorithm. It means that RH-ADMM renders the solution of the D-NLP problem fully decentralized for which modern CPU architectures can be utilized.  In general, based on the communication topology of the RH-ADMM algorithm, it is convenient to think of RH-ADMM as a message-passing algorithm on a hypergraph, where each node corresponds to a subsystem and hyperedges correspond to shared variables. As mentioned, a well-known protocol for implementing this form of communication structure for parallel algorithms is the MPI. In particular, the MPI  is a language-independent message-passing specification used for parallel algorithms, which is the most widely used model for high-performance parallel computing nowadays. There are numerous implementations of MPI on a variety of distributed platforms, and interfaces to MPI are available from a wide variety of languages, including {\tt C/C++}, {\tt Java}, and {\tt Python}.

There are multiple strategies to implement RH-ADMM in MPI. One simple way to do so is by coding the algorithm using the \textit{Single Program, Multiple Data} (SPMD) programming style, in which each processor or subsystem runs the same program code, but each has its own set of local variables, objective functions, and data and operates in a separate subset of the data. 
\begin{algorithm}
\begin{small}
\SetAlgoLined
\KwResult{Initial local binary vector $\zz_j, j = 1,...,K$ }
\caption{SFP algorithm}
1. solve problem \eqref{dist-prim} in parallel and obtain $\x_i^0,\, \forall i = 1,...,N$, and $\y_j^0,\, \forall j = 1,...,K$;\\
2. $\Bar{\y}_j^0\xleftarrow{} \Pi_{C_j}(\y_j^0)$, $\forall j = 1,...,K$\\
3. generate $\zz_j^0\in\{0,1\}^n, \, \forall j = 1,...,K$, based on feasible vector  $\bar{\mathbf{y}}_j^0$;\\
4. $reduce$ $f_i(\x_i^0)$, \ 
$ub^0\xleftarrow{} \sum\limits_{i = 1}^{N}f_i(\x_i^0)$;\\
\label{sfp:alg}
\end{small}
\end{algorithm}

\subsection{DiPOA Algorithm}
\begin{algorithm}
\begin{small}
\SetAlgoLined
\KwResult{The solution of the SCP problem \eqref{dis-ccp}}
\caption{DiPOA algorithm}
$SFP \xleftarrow{} True/False$ (whether or not SFP is used);\\
$SOC \xleftarrow{} False$;\\
\eIf{SFP}{Run Algorithm \ref{sfp:alg} to  obtain initial local binary vectors $\zz_j, \, \forall j = 1,...,K$;
}{Initialize the lower bound, $lb^k$ and the upper bounds, $ub^k$, and generate arbitrary feasible binary vectors $\zz_j^k,\, \forall j = 1,...,K$;\\
}
compute $r^0$ from \eqref{rgap};\\
$k \xleftarrow{} 0$;\\
\While{$ r^k \geq \epsilon$}{
$k\xleftarrow{} k + 1$;\\
1. solve problem D-NLP  \eqref{dist-prim} in parallel for fixed $\zz_j^{k-1},\, \forall j = 1,...,K$, using the RH-ADMM Algorithm; \\
2. $reduce$ $f_i(\x_i^0)$ and $gather$ local information to the central cut storage unit.\\
\eIf{SOC}{
 3. generate \textbf{second order underestimators} \eqref{soc} around $\x_i^k,\, \forall i = 1,...,N$, in parallel;\\
}{
3. generate \textbf{linear underestimators}, \eqref{lin-app-feasb}, around $\x_i^k,\, \forall i = 1,...,N$,  in parallel;\\
}
4. solve the master problem \eqref{relaxed-master-miqcp} $bcast$ \ $\zz_j^k$ \\
5. Update $lb^k$, $ub^k$, $r^k$ \eqref{rgap}, and $e^k$ \eqref{event}\\
\eIf{$e^k \leq tol$}{
6. $SOC \xleftarrow{} True$;
}{6. $SOC \xleftarrow{} False$;
}
}
\label{alg1}
\end{small}
\end{algorithm}
We provide two algorithms to perform the computations. Algorithm \ref{sfp:alg} implements the specialized feasibility pump, while Algorithm \ref{alg1} corresponds to the DiPOA methodology. 

Algorithm \ref{sfp:alg} consists of four main steps. In step $1$,  problem \eqref{dist-prim-sfp} is solved and then its solution is projected onto the sparsity constraint set, $\mathcal{C}_j$. We denote $\Bar{\y}_j$ as the projected vectors. Once the projection step is performed,  local binary vectors $\zz_j$ are obtained, $j = \{1,...,K\}$, as it was illustrated in Example 1.  The computations of this step are performed by the primal level (as was shown in Figure \ref{num-hier}), which solves problem \eqref{dist-prim-sfp} cooperatively until consensus is achieved. Then the optimal local solutions are sent to the LFC level, where the local projections are performed and local binary vectors are generated. In the fourth step, the \textit{reduce} operation is performed, where the reduce refers to {\tt MPI\_Reduce}. In this step, the sum of all local objective functions is computed and sent to the {\tt Root} node, where the master's problem is solved. Lastly, using the results of this step, we initialize the upper bound.

Algorithm \ref{alg1} starts with an initialization step. In this phase, we initialize the DiPOA algorithm by choosing whether or not SFP and/or SOC are used. We first initialize the algorithm without the SoCut property, since the event-triggered scheme is utilized to switch between linear and second-order cuts. In the case that SFP is selected, Algorithm \ref{sfp:alg} is used to generate an initial feasible binary vector $\zz_j$. Otherwise, we generate arbitrary feasible binary vectors $\zz_j$. By using the initial upper bound and lower bound, we obtain the initial relative duality gap according to equation \eqref{rgap}. After the initialization phase is terminated, we start iterating the algorithm by checking the relative optimality gap computed at each iteration $k$, which is induced by the upper and lower bound. We then update the upper bound by solving the distributed NLP problem \eqref{dist-prim} for fixed local binary vectors $\zz_j$. Similar to the SFP step, these computations are performed by the RH-ADMM algorithm and in parallel. As before, the primal level is responsible for calculating the upper bound and also the approximation points provided by the primal problem \eqref{dist-prim}. Once a consensus is reached, the reduce operation is performed and the upper bound is updated. Step $4$ consists of sending the local information (\textit{i.e.}, local solutions, gradient, and minimum eigenvalue) to the cut manager where this information is processed. The \textit{gather} operation denotes 
{\tt MPI\_Gather}  which is an MPI routine to transmit the local data to the {\tt  Root}. Step $5$ is responsible for the \textit{synchronization}. This step guarantees that all local nodes compute and transmit their local information to the {\tt Root}. Here  \textit{Barrier} refers to the {\tt MPI\_Barrier}  routine which performs node synchronization. 

After generating cuts according to their order and storing their information, the master's problem \eqref{relaxed-master-miqcp} is solved with a MIP solver. The cut generation procedure is performed by the cut manager. Note that the number of cuts generated by a particular cut manager is the same as the number of nodes in the corresponding LFC unit. For instance, if LFC$_j$ consists of $N_j$ nodes, then CG$_j$ generates $N_j$ cuts for the master's problem. Since we initialize DiPOA without the SoCut property, at the first iteration, a linear cut is generated and added to the master's problem. Then, problem \eqref{relaxed-master-milp} is solved to achieve new binary vectors and also the lower bound on the objective value. Once the master's problem is solved, the {\tt Root}  node \textit{Broadcasts} the obtained binary vector $\zz_j^k$ to the local LFCs, through the {\tt MPI\_Bcast}  routine. By using the updated information about the upper and lower bound, we update the relative gap at the current iteration and compute the event, $e^k$ from \eqref{event}. If $e^k$ is smaller than the tolerance $tol$, we use SoCuts. Otherwise, we continue with linear cuts. This procedure is repeated until the relative duality gap is smaller than a given tolerance of $\epsilon$.  

In addition to the algorithmic details that were explained above, some other numerical aspects must be taken into account, which are discussed below.

\subsection{Numerical Considerations on Linear and Quadratic Cuts}
In this section, we discuss some numerical aspects of the cut generation procedure, which is performed by the cut manager of the numerical hierarchy depicted in Figure \ref{num-hier}. To produce the linear cuts, first, the cut manager receives the cut information from the LFC level (after RH-ADMM converged) and then, stores them in the cut-storage pool. In particular, each local cut manager unit $j$ receives and stores the local solution $\x_i^k,\, \forall i \in \mathcal{E}_j$, the optimal value $f_i(\x_i^k), \, \forall i \in \mathcal{E}_j$, and the local gradient  $\nabla f_i(\x_i^k), \, \forall i \in \mathcal{E}_j$. To add a cut to the master's problem, the cut information is retrieved from the cut-storage pool and added to the master problem. Otherwise, the information is kept inside the storage pool. For the SoCuts', the cut information also contains the strong convexity parameter $m_i$. Considering one of the cases in which $m_i$ can be computed efficiently, there is no need to compute and store the Hessian matrix at each iteration $k$ of the DiPOA algorithm. More importantly, storing the Hessian matrices is not efficient from the MPI point of view. For instance, in large-scale instances of problem \eqref{dis-ccp}, broadcasting the Hessian matrix at each iteration stands as a bottleneck for the entire algorithm. Hence, it is essential to reduce the amount of information that is broadcast amongst the CNs. Therefore, from the perspective of distributed numerical computation, the SoCuts \eqref{soc} can provide better performance. The only limitation is, however, that calculating the strong convexity parameter, $m$, for the nonlinear functions is not straightforward. As mentioned before, for most of the learning and control problems such as classification, optimal and predictive control, and regression, $m$ can be computed efficiently and accurately. 

\section{Computational Experiments and Algorithm Evaluation}\label{exper}
This section evaluates the accuracy and the performance of DiPOA by applying Algorithm \ref{alg1} to solve two SCP problem instances, namely DSLR and SQCQP problems. Results from a comparison with state-of-the-art MINLP solvers are reported to validate the optimality of DiPOA.
Further, several instances of these problems are solved, for different settings and scenarios, to analyze the efficiency of the DiPOA algorithm.
Since several centralized and distributed applications of SCP problems are found in applied statistics, machine learning, and control areas, in this paper we consider that DiPOA will be mostly applied for solving sparse convex programming problems arising in these fields, although it is possible to apply it in more generalized scenarios that fit the SCP problem  \eqref{D-CCP}. Therefore, in this section, we focus on sparse classification and model fitting problems.

\subsection{Implementation Details}
The experiments are performed on a machine consisting of two {\tt Xeon(R)} {\tt CPU}  {\tt E5-2630}  {\tt 2.20GH}  processors, each of which with $10$ logical cores and two-way hyper-threaded, and $64$ GB of memory. The DiPOA algorithm is implemented in {\tt Python}   programming language and the linear algebra operations are performed using {\tt Numpy},  which is a scientific computing library of {\tt Python}. Moreover, the MPI and the message-passing paradigm are implemented by {\tt Python}  {\tt mpi4py}  library which is available for free. Finally, the  MILP master's problem of the DIPOA algorithm is implemented and solved with {\tt Gurobi}  {\tt Python}  API \cite{gurobi}. The source code of the algorithm is also available on \textcolor{black}{Github}\footnote{\url{https://github.com/Alirezalm/dccp.git}.}.

\subsection{Distributed Logistic Regression Problem}
Sparse classification is one of the central problems in statistics and machine learning as it leads to more interpretable models. One of the most popular classification methods in machine learning is \textit{logistic regression} where the response vector has only two states.
In machine learning problems, when the number of potential features, $n$, is large, it is often desired to identify a critical subset of features, which are primarily responsible and contribute to the response. This leads to training more \textit{interpretable} models and improves \textit{prediction accuracy} by eliminating unnecessary variables and increasing the model's generalizability. For this reason, the sparsity requirement is incorporated into logistic regression models. In principle, this approach involves identifying a subset of $\kappa < n$ features that are related to the response. Then, a logistic regression model is fitted on the reduced set of variables. This problem is also known as \textit{subset selection} or Sparse Logistic Regression (SLR) problem \cite{james2013introduction}. This requirement can be formulated by a sparsity constraint framework. Hence, the centralized SLR problem is formulated as the optimization problem given in Eq. \eqref{SLR}.
First, we focus on the centralized and dense formulation of the SLR problem, and then we move to the distributed set-up while introducing the sparsity property in the model. Given a response vector $\Gam \in \{-1,1\}^p$, a training dataset matrix $\X \in \mathcal{R}^{p \times n}$, and unknown regression coefficients $\thh \in \mathcal{R}^n$, the aim is to build a linear classification model from the training samples by maximizing the probability of each class. For the centralized case, the SLR problem consists of \textcolor{black}{solving the following SCP problem:}
\begin{subequations}\label{SLR}
\begin{align}
    \min_{\thh} &\sum_{\ell = 0}^p \log\left[1 + \mathrm{e}^{-(\x_{\ell}^T\thh)\Gamma_{\ell}}\right] + \frac{\lambda}{2}\norm{\thh}_2^2\\
s.t.~    & \norm{\thh}_0 \leq \kappa \label{cardlr}
\end{align}
\end{subequations}
where $\x_{\ell}$ is the $\ell$-th row of $\X$, $\Gamma_{\ell}$ it the $\ell$-th element of $\Gam$, and $\lambda >0$ is the regularization parameter. Without the sparsity constraint, problem \eqref{SLR}  is an unconstrained centralized convex optimization problem, whose optimal solution can be found using unconstrained optimization algorithms, such as gradient descent methods (\textit{e.g.}, batch gradient descent and stochastic gradient descent) and Newton's method. 
By denoting $\Gam_i \in \{-1,1\}^p$ as  the $i$-th local response vector and $\X_i$ as the $i$-th local training dataset matrix, the distributed formulation of Problem \eqref{SLR}, namely the Distributed Sparse Logistic Regression (DSLR) problem, for $N$ local nodes, can be written as 
\begin{subequations}\label{DSLR}
\begin{align}
      \min_{\thh} &\sum_{i = 1}^{N} \sum_{\ell = 0}^p \log\left[1 + \mathrm{e}^{-(\x_{i,\ell}^T\thh)\Gam_{i,\ell}}\right] + \frac{\lambda}{2}\norm{\thh}_2^2\\
  s.t.~ & \norm{\thh}_0 \leq \kappa
\end{align}
\end{subequations}
where $N$ is the number of different locations where the data is stored,  $\x_{i,\ell}$ is the $\ell$-th row of $\X_i$, and $\Gam_{i,\ell}$ is the $\ell$-th element of $\Gam_i$. 
\begin{remark}
The use of ridge regularization is justified for multiple reasons. It enforces the uniqueness of the fitted coefficients and therefore accounts for over-determined systems or multi-collinearities. Otherwise, the coefficients would gain several degrees of freedom, which could result in the coefficients having high variance. Another reason for the use of ridge penalization is that a model selection process can generate selection bias which leads to an overestimation of the coefficients. The ridge regression causes coefficients to shrink and can therefore correct the additional bias. Moreover, from the numerical point of view, the ridge parameter ensures the strong convexity of the log-likelihood function and thus improves the condition number of the Hessian of the objective function, which can help the first-order algorithms inside RH-ADMM to converge faster.
\end{remark}
Following the DiPOA steps, we reformulate problem \eqref{DSLR} as the following D-MINLP problem:
\begin{subequations}\label{DSLR-admm}
\begin{align}
      \min_{\thh_i, \y_j, \zz_j} &\sum_{i = 1}^{N} \sum_{\ell = 0}^p \log\left[1 + \mathrm{e}^{-(\x_{i,\ell}^T\thh_i)\Gam_{i,\ell}}\right] + \frac{\lambda}{2}\norm{\thh_i}_2^2\\
      & \thh_i = \y_j, \, \forall i \in \mathcal{E}_j, \, \mathcal{E}_j \in \mathcal{E}, \, \forall j=1,...,K\\
          &-M_j\zz_j \leq \y_j \leq M_j \zz_j\\
    &e_j^T\zz_j \leq \kappa
\end{align}
\end{subequations}
Problem \eqref{DSLR-admm} belongs to a subclass of the SCP problem \eqref{dis-ccp} to which the DiPOA algorithm can be applied. In the following section, we provide the results obtained by solving problem \eqref{DSLR-admm} with Algorithm \ref{alg1}.
\subsubsection*{Numerical Results: }
To evaluate DiPOA for the DSLR problem, we generate (or split the original dataset into) $N$ random local datasets which are standardized to have zero
mean and unit $\ell_2$ norm. 
 The response vector $\Gamma$ is generated according to the logistic function as follows,
\begin{equation*}
    \Gam_{i,\ell} = \text{round}\left (\frac{1}{1+\text{exp}(-\theta^T\x_{i,\ell})}\right ), \, \forall \ell=0,\dots,p,
\end{equation*} Each local dataset has the same number of features compared to other local datasets, however, the number of samples is reduced. Taking into account that most MINLP solvers are not distributed, we compare DIPOA with the state-of-the-art MINLP centralized solvers. Accordingly, we consider BONMIN, KNITRO, DICOPT, and SHOT through GAMS interface.  We also compared the DIPOA evaluation results to the OA algorithm proposed by \cite{Bertsimas2017}, which was proposed for centralized sparse logistic regression problems.

\begin{figure}[htb!]
    \centering
    \includegraphics[width=0.55\textwidth]{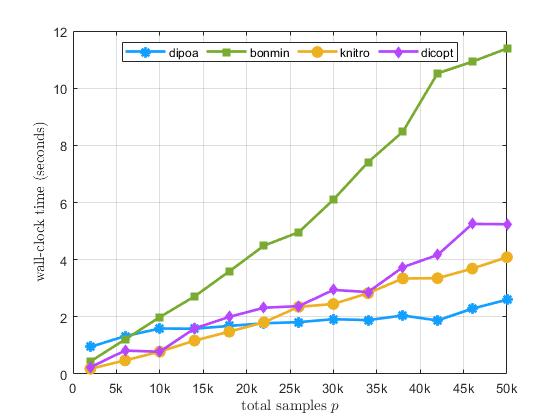}
    \caption{\footnotesize Comparison between DiPOA, BONMIN, KNITRO, and DICOPT} 
    \label{dslr:comparison}
\end{figure}

We evaluate the DiPOA based on five main scenarios with different problem settings and parameters. Each scenario consists of solving multiple instances of the DSLR problem. In the first scenario, \textbf{SC-I}, we consider a different number of total sample points and evaluate DiPOA for the case where only SFP is available and also for the case where both SFP and SoCut features are activated. In the second scenario, \textbf{SC-II}, we fix the total number of sample points and change the number of decision variables. Finally, in the last three scenarios, \textbf{SC-III - SC-V}, we compare fully-featured DiPOA with MINLP solvers for relatively small, medium, and large problem instances. Each scenario consists of different settings which are provided in Table \ref{sc:settings}. In this table, {\tt total-samples} refers to the number of data points, {\tt num-var} and {\tt num-nodes} denote the number of variables and nodes, respectively.
\begin{table}[htb!]
    \centering
    \caption{\footnotesize Scenario settings \label{sc:settings}}
    {    \begin{tabular}{SSSSS}
    \toprule
            {\Verb!scenario!} &{ \Verb!total-samples!} & { \Verb!num-var!} & { \Verb!acceptable-sparsity [$\%$]!} & { \Verb!num-nodes!}\\\toprule
            {\textbf{SC-I}} &{$[2k \mathrel{{.}\,{.}} 50k]\nobreak$} & {$20$} & {$25$} & {$10$}\\
            {\textbf{SC-II}} &{$50k$} & {$40-200$} & {$25$} & {$10$}\\\midrule
            {\textbf{SC-III}} &{$10k$} & {$20$} & {$[10\mathrel{{.}\,{.}}\nobreak 90]$} & {$10$}\\
            {\textbf{SC-IV}} &{$100k$} & {$200$} & {$[10 \mathrel{{.}\,{.}}\nobreak 90]$} & {$10$}\\
            {\textbf{SC-V}} &{$300k$} & {$300$} & {$[10 \mathrel{{.}\,{.}}\nobreak 90]$} & {$10$}\\
    \bottomrule
    \end{tabular}
    }
    {}
\end{table}
The numerical results from \textbf{SC-I} appear in Table \ref{sc-i:res}. 
 In this case, DiPOA is evaluated according to two different settings. In the first setting (\Verb!dipoa-sfp! column), DiPOA only generates linear cuts, and only SFP is activated, whereas, in the second setting (\Verb!dipoa-soc! column),  SoCuts can be generated. In this Table and subsequent tables, \Verb!rel-gap! refers to relative optimality gap \eqref{rgap} and \Verb!time! is CPU wall-clock time. It can be seen in  Table \ref{sc-i:res} that second-order cuts can drastically improve the algorithm convergence in terms of wall-clock time and the number of cuts. Moreover, we observe that as the number of all samples increases, the number of FOC cuts and therefore the wall-clock time decreases. In other words, a large number of sample points often lead to a smaller number of iterations, however, with more computational complexity per iteration. Considering the centralized architecture, although the number of iterations can be small, the computational burden is larger than the distributed case since the entire dataset is used in the evaluation of the objective function. This behavior can be seen in Figure \ref{dslr:comparison} where the scalability of the DiPOA is evaluated. This figure compares the wall-clock time between DiPOA and BONMIN, KNITRO, and DICOPT which are centralized MINLP solvers. Compared to the centralized solvers, for the problem instances with a relatively small number of data points ($p \leq 6k$), DiPOA needs more time to provide the solution. In contrast, as we increase the size of the dataset,  DiPOA behaves more robustly and needs a smaller execution time. Hence, DiPOA can scale well for problem instances with much larger data sets. We also note that SHOT was able to successfully solve only the first instances ($2k$, $6k$, and $10k$) of problem SC-I,  while for larger instances and other scenarios it reaches the cut-off time of $600$ seconds. For the first and second instances (where numbers of total samples are $ 2k $ and $6k$) the solution time is $57.86$ and $286.63$ seconds. For the third case with $10k$ as the number of total samples, SHOT converged in $536.10$ seconds. Since SHOT failed to solve the larger problem instances within the considered cut-off time, we omit its results in the figures and tables. 
 \begin{table}[htb!]
    \centering
    \caption{\footnotesize Numerical Results for SC-I (Varying Number of Samples) \label{sc-i:res}}
    {    \begin{tabular}{c SSSSS SSSSS}
    \toprule
    \multirow{1}{*}{} &
     \multicolumn{5}{c}{\Verb!dipoa-sfp!} &
     \multicolumn{5}{c}{\Verb!dipoa-soc!} \\
     \cmidrule(lr){2-6} \cmidrule(lr){7-11}
     {\Verb!samples! }& {\Verb!time!}&{\Verb!nfoc!}&{\Verb!gap! }&{\Verb!mip! }&{\Verb!nlp!}& {\Verb!time!}&{\Verb!nsoc!}&{\Verb!gap! }&{\Verb!mip! }&{\Verb!nlp! }\\
     \toprule
            {$2k$}    &{$600$}      &{  $1510$ } & {$0.133$ }  &{$ 439.39$}&{$160.60 $}& { $0.947$}& { $20$}& {$0.063$ }&{$0.79 $}&{$ 0.14$}\\
            {$6k$}    &{$600 $}     &{  $780 $ } & {$ 0.0123$ }&{$439.39 $}&{$160.60 $}& { $1.326 $}& { $10 $}& {$0.38 $ }&{$ 1.10$}&{$ 0.21$}\\
            {$10k$}   &{$336.141 $} &{  $ 430$ } & {$ 0.896$ }&{$ 246.16$}&{$89.97 $} & { $ 1.591$}& { $10 $}& {$0.36 $ }&{$ 1.32$}&{$ 0.26$}\\
            {$14k$}   &{$ 268.325$} &{  $ 361$ } & {$0.899 $ } &{$196.49 $}&{$71.82 $}& { $1.585 $}& { $10 $}& {$0.13 $ }&{$ 1.32$}&{$ 0.26$}\\
            {$18k$}   &{$ 186.343$} &{  $270 $ } & {$0.864 $ }&{$136.464 $}&{$49.87 $}& { $ 1.681$}& { $10 $}& {$0.13 $ }&{$ 1.40$}&{$ 0.27$}\\
            {$22k$}   &{$152.284 $} &{  $200 $ } & {$ 0.896$ }&{$111.52 $}&{$40.761 $} & { $1.775 $}& { $10 $}& {$0.067 $ }&{$1.48 $}&{$ 0.29$}\\
            {$26k$}   &{$147.808 $} &{  $180 $ } & {$ 0.859$ } &{$108.23 $}&{$39.56 $}& { $ 1.814$}& { $10 $}& {$0.074 $ }&{$1.51 $}&{$ 0.29$}\\
            {$30k$}   &{$122.878 $} &{  $170 $ } & {$0.867 $ } &{$88.75 $}&{$34.12 $}& { $1.918 $}& { $10 $}& {$0.042 $ }&{$1.60 $}&{$ 0.31$}\\
            {$34k$}   &{$120.662 $} &{  $ 160$ } & {$0.889 $ } &{$ 87.15$}&{$ 33.50$}& { $1.886 $}& { $10 $}& {$0.035 $ }&{$ 1.57$}&{$ 0.31$}\\
            {$38k$}   &{$ 118.914$} &{  $150$ }  & {$0.897 $ } &{$ 85.88$}&{$ 33.02$}& { $2.046 $}& { $10 $}& {$0.073 $ }&{$ 1.70$}&{$0.33 $}\\
            {$42k$}   &{$ 86.67$}   &{  $ 110$ } & {$0.885 $ } &{$ 62.60$}&{$ 24.06$}& { $1.874 $}& { $10 $}& {$0.067 $ }&{$ 1.56$}&{$ 0.30$}\\
            {$46k$}   &{$ 11.041$}  &{  $30 $ }  & {$0.708 $ } &{$11.04 $}&{$ 3.06$}& { $2.286 $}& { $10 $}& {$0.051 $ }&{$1.90 $}&{$ 0.37$}\\
            {$50k$}   &{$11.413$}   &{  $30 $ }  & {$0.806 $ } &{$8.24 $}&{$ 3.16$}& { $2.598 $}& { $10 $}& {$0.082 $ }&{$ 2.17$}&{$ 0.42$}\\
    \bottomrule
    \end{tabular}
    }
    {}
\end{table} 

\begin{figure}[htb!]
    \centering
    \includegraphics[width=0.55\textwidth]{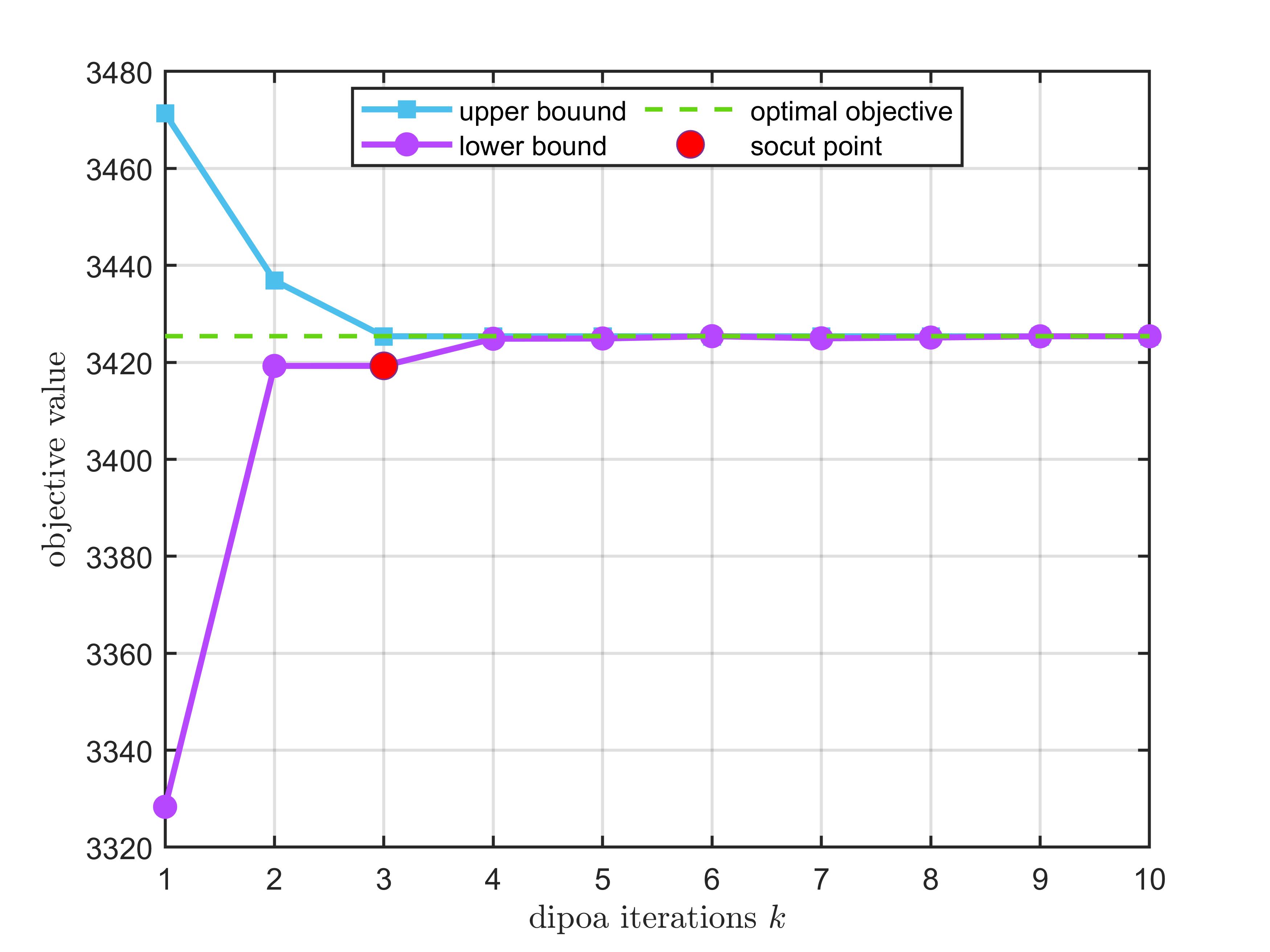}
    \caption{\footnotesize The impact of the SoCuts and the event-triggered scheme.} 
    \label{ublb}
\end{figure}

 Table \ref{sc-ii:res} provides the numerical results for SC-II. In this case, the number of data points is fixed to be $50k$ and DiPOA is evaluated concerning a varying number of variables, which runs with both SFP and SoCuts activated. These results show that DiPOA is robust to the increase of variables. 
 \begin{table}[htb!]
    \centering
    \caption{\footnotesize Numerical Results for SC-II (Varying Number of Variables) \label{sc-ii:res}}
    {    \begin{tabular}{c SS SS SS SS SS}
    \toprule
    \multirow{2}{*}{} &
     \multicolumn{2}{c}{\Verb!dipoa!} &
     \multicolumn{2}{c}{\Verb!bonmin!}&
     \multicolumn{2}{c}{\Verb!knitro!}&
      \multicolumn{2}{c}{\Verb!dicopt!}\\
     \cmidrule(lr){2-3} \cmidrule(lr){4-5} \cmidrule(lr){6-7} \cmidrule(lr){8-9}
     {\Verb!num-var! }& {\Verb!time[sec]!}& {\Verb!rel-gap!$[\%$]}& {\Verb!time[sec]!}& {\Verb!rel-gap$[\%$]!}& {\Verb!time[sec]!}& {\Verb!rel-gap!$[\%$]} & {\Verb!time[sec]!}& {\Verb!rel-gap!$[\%$]}\\
     \toprule
            {$40$}  &{  $2.446 $ } & {$0.06 $ }  & { $25.561 $}& { $0.12 $}& {$ 11.637$ }&{ $0.012 $}& { $17.188 $}& {$0.22 $ } \\
            {$60$}  &{  $2.651 $ } & {$0.051 $ }  & { $59.001 $}& { $0.09 $}& {$20.335 $ }&{ $0.013 $}& { $22.900 $}& {$0.087 $ }  \\
            {$80$}  &{  $13.872 $ } & {$0.003 $ }  & { $97.167 $}& { $0.13 $}& {$30.336 $ }&{ $0.063 $}& { $32.403 $}& {$0.31 $ }  \\
            {$100$} &{  $14.321 $ } & {$0.005 $ }  & { $136.693 $}& { $0.18 $}& {$40.33 $ }&{ $0.036 $}& { $55.151 $}& {$0.02 $ }  \\
            {$120$} &{  $15.005 $ } & {$0.005 $ }  & { $203.814 $}& { $0.02 $}& {$49.354 $ }&{ $0.063 $}& { $60.003 $}& {$ 0.03$ } \\
            {$140$} &{  $15.537 $ } & {$ 0.009$ }  & { $272.458 $}& { $0.19 $}& {$56.332 $ }&{ $0.055 $}& { $114.103 $}& {$0.06 $ }  \\
            {$160$} &{  $16.315 $ } & {$ 0.006$ }  & { $336.524 $}& { $0.15 $}& {$67.356 $ }&{ $0.041 $}& { $117.705 $}& {$0.08 $ }  \\
            {$180$} &{  $17.189 $ } & {$0.003 $ }  & { $390.053 $}& { $0.36 $}& {$76.387 $ }&{ $0.032 $}& { $164.234 $}& {$0.03 $ }  \\
            {$200$} &{  $18.002 $ } & {$0.007 $ }  & { $463.526 $}& { $0.45 $}& {$85.569 $ }&{ $0.045 $}& { $228.631 $}& {$0.08 $ }  \\
    \bottomrule
    \end{tabular}
    }
    {}
\end{table}
The numerical results of SC-III are presented in Table \ref{sc-iii:res}. 
  The experiment aims to analyze the sensitivity of DiPOA for the number of non-zero elements in the solution (bound $\kappa$). Since the problem size is small, the centralized solvers outperform DiPOA in terms of execution time. 
   \begin{table}[htb!]
    \centering
    \caption{\footnotesize Numerical Results for SC-III (Varying Number of Non-Zeros) \label{sc-iii:res}}
    {    \begin{tabular}{c SS SS SS SS}
    \toprule
    \multirow{2}{*}{} &
     \multicolumn{2}{c}{\Verb!dipoa!} &
     \multicolumn{2}{c}{\Verb!bonmin!}&
     \multicolumn{2}{c}{\Verb!knitro!}&
     \multicolumn{2}{c}{\Verb!dicopt!}\\
     \cmidrule(lr){2-3} \cmidrule(lr){4-5} \cmidrule(lr){6-7} \cmidrule(lr){8-9}
     {$\kappa$ }& {\Verb!time[sec]!}& {\Verb!rel-gap!}& {\Verb!time[sec]!}& {\Verb!rel-gap!}& {\Verb!time [sec]!}& {\Verb!rel-gap!} & {\Verb!time [sec]!}& {\Verb!rel-gap!}\\
     \toprule
            {$ 2$}  &{  $2.927 $ } & {$0.075 $ }  & { $2.115 $}& { $0.15 $}& {$0.861 $ }&{ $0.12 $}& { $1.095 $}& {$ 0.12$ }\\
            {$4 $}  &{  $2.746 $ } & {$0.120 $ }  & { $1.970 $}& { $0.061 $}& {$ 0.860$ }&{ $0.204 $}& { $1.119 $}& {$0.20 $ }\\
            {$6 $}  &{  $2.862 $ } & {$ 0.160$ }  & { $1.950 $}& { $ 0.102$}& {$0.866 $ }&{ $0.21 $}& { $1.192 $}& {$ 0.16$ }\\
            {$8$}  &{  $3.164 $ } & {$0.150 $ }  & { $2.021 $}& { $ 0.23$}& {$0.783 $ }&{ $ 0.13$}& { $1.365 $}& {$0.12 $ }\\
            {$ 12$}  &{  $2.816 $ } & {$0.249 $ }  & { $2.024 $}& { $0.12 $}& {$0.673 $ }&{ $0.132 $}& { $1.157 $}& {$0.18 $ }\\
            {$14 $}  &{  $2.936 $ } & {$0.200 $ }  & { $2.040 $}& { $0.10 $}& {$0.715 $ }&{ $0.153 $}& { $1.490 $}& {$0.16 $ }\\
            {$16 $}  &{  $2.729 $ } & {$0.250 $ }  & { $2.034 $}& { $0.17 $}& {$ 0.656$ }&{ $0.106 $}& { $1.021 $}& {$0.13 $ }\\
            {$18 $}  &{  $2.459 $ } & {$0.173 $ }  & { $2.094 $}& { $0.09 $}& {$ 0.686$ }&{ $0.151 $}& { $1.156 $}& {$0.15 $ }\\
    \bottomrule
    \end{tabular}
    }
    {}
\end{table}
Finally, Table \ref{sc:iv-v} presents the numerical results for medium and large scale scenarios (\textit{i.e.}, SC-IV and SC-V). In these scenarios, the centralized solvers were not capable of providing a solution for the problem instances. The properties of the problem instances are defined in Table \ref{sc:settings} for each scenario.
\begin{table}[htb!]
    
    \caption{\footnotesize Numerical Results for SC-IV and SC-V (Medium and Large Scenarios)}\label{sc:iv-v}
    \begin{subtable}[t]{0.4\textwidth}
    \centering
    \caption{\footnotesize \textbf{SC-IV} \label{sc-iv:res}}
    {    \begin{tabular}{SSS}
    \toprule
            {$\kappa$} &{ \Verb!dipoa-time[sec]!} & { \Verb!dipoa-gap[$\%$]!}\\\toprule
            {$20 $} &{ $6.608 $ } & {$ 0.11$} \\
            {$40 $} &{ $ 6.787$ } & {$0.16 $} \\
            {$60 $} &{ $6.945 $ } & {$0.14 $} \\
            {$80 $} &{ $6.550 $ } & {$0.17 $} \\
            {$100 $} &{ $7.078 $ } & {$0.20 $} \\
            {$120 $} &{ $ 6.595$ } & {$ 0.10$} \\
            {$160 $} &{ $6.601 $ } & {$ 0.11$} \\
            {$180 $} &{ $5.586 $ } & {$ 0.22$} \\
    \bottomrule
    \end{tabular}
    }
    {}
    \end{subtable}
    \hfill
    \begin{subtable}[t]{0.4\textwidth}
    \centering
    \caption{\footnotesize \textbf{SC-V} \label{sc-v:res}}
    {    \begin{tabular}{SSS}
    \toprule
            {$\kappa$} &{ \Verb!dipoa-time[sec]!} & { \Verb!dipoa-gap[$\%$]]!}\\\toprule
            {$ 20$} &{ $18.699 $ } & {$0.03 $} \\
            {$60 $} &{ $18.424 $ } & {$0.05 $} \\
            {$ 100$} &{ $18.598 $ } & {$0.06 $} \\
            {$140 $} &{ $18.402 $ } & {$ 0.11$} \\
            {$180 $} &{ $17.618 $ } & {$0.092 $} \\
            {$220 $} &{ $18.785 $ } & {$0.097$} \\
            {$ 240$} &{ $18.582 $ } & {$ 0.095$} \\
            {$ 260$} &{ $18.201 $ } & {$0.093 $} \\
    \bottomrule
    \end{tabular}
    }
    {}
    \end{subtable}

\end{table}
Additionally, we investigate the impact of SoCuts and their event-triggered scheme on the solution of the master's problem.  Figure \ref{ublb} illustrates the impact on the upper and the lower bounds in the context of the DSLR problem.
In particular, we allow DiPOA to perform $10$ iterations for a small instance of the DSLR problem with $\kappa = 5$, $\theta = 10$, and $p = 2000$. In Figure \ref{ublb}, purple circles indicate the lower bound and the black squares show the upper bound generated by DiPOA. The red circle is the point at which quadratic cuts are generated. As the figure depicts, the  SoCuts improve the convergence of the algorithm by inducing a higher lower bound. Moreover, the event-triggered scheme prevents the algorithm from generating a high number of  SoCuts. For example, in this problem instance, only one SoCut is generated, preventing the master's problem from becoming overly complex with the quadratic constraints. Although this procedure provides a better lower bound while controlling the complexity of the master's problem, it remains an unknown function of the problem data. 
Finally, numerical results comparing DiPOA to the Centralized Outer Approximation (COA) of \cite{Bertsimas2017} are provided in Table \ref{comparison-dslrr}. The same settings from \cite{Bertsimas2017} are adopted. As the table shows, the DiPOA algorithm outperforms the centralized architecture and solves the problem instances in significantly less computational time. 
 \begin{table}[htb!]
    \centering
    \caption{\footnotesize Numerical Results Comparing DiPOA with Centralized OA for $\kappa = 5$. \label{comparison-dslrr}}
    {    \begin{tabular}{SSSSS}
    \toprule
            {\Verb!total-samples!} &{ \Verb!num-var!} & { \Verb!num-nodes!} & { \Verb!DiPOA-time!\,[\textcolor{black}{sec}]} & { \Verb!COA-time!\,[sec]}\\\toprule
            {$100$}  & {$10$} & {$4$} &\textcolor{black}{{$0.077$}} & {$<1$} \\  
        {$1k$}  & {$100$} & {$4$} &\textcolor{black}{{$2.581$}} & {$15$} \\  
        {$2k$}  & {$200$} & {$4$} &\textcolor{black}{{$3.298$}}  & {$16$} \\
    \bottomrule
    \end{tabular}
    }
    {}
\end{table}

\subsection{Sparse Quadratically Constrained Quadratic Programming Problem}
In this section, we provide numerical evaluation results for the SQCQP problem, which is a SCP problem with quadratic objective and constraint functions,  namely:
\begin{subequations}\label{dis-dsqcp}
\begin{align}
 \min_{\x \in \mathcal{R}^n}\,\,  &\sum_{i=1}^{N}\left (\frac{1}{2}\x^TQ_i\x + q_i^T\x + d_i \right )\\
  s.t.~  &\frac{1}{2}\x^TP_h\x + c_h^T\x  + r_h\leq 0, \,\, \forall h=1,...,m \\
    & \x \in \Omega \\
    &\norm{\x}_0 \leq \kappa 
\end{align}
\end{subequations}

\begin{figure}
    \centering
    \includegraphics[width=0.55\textwidth]{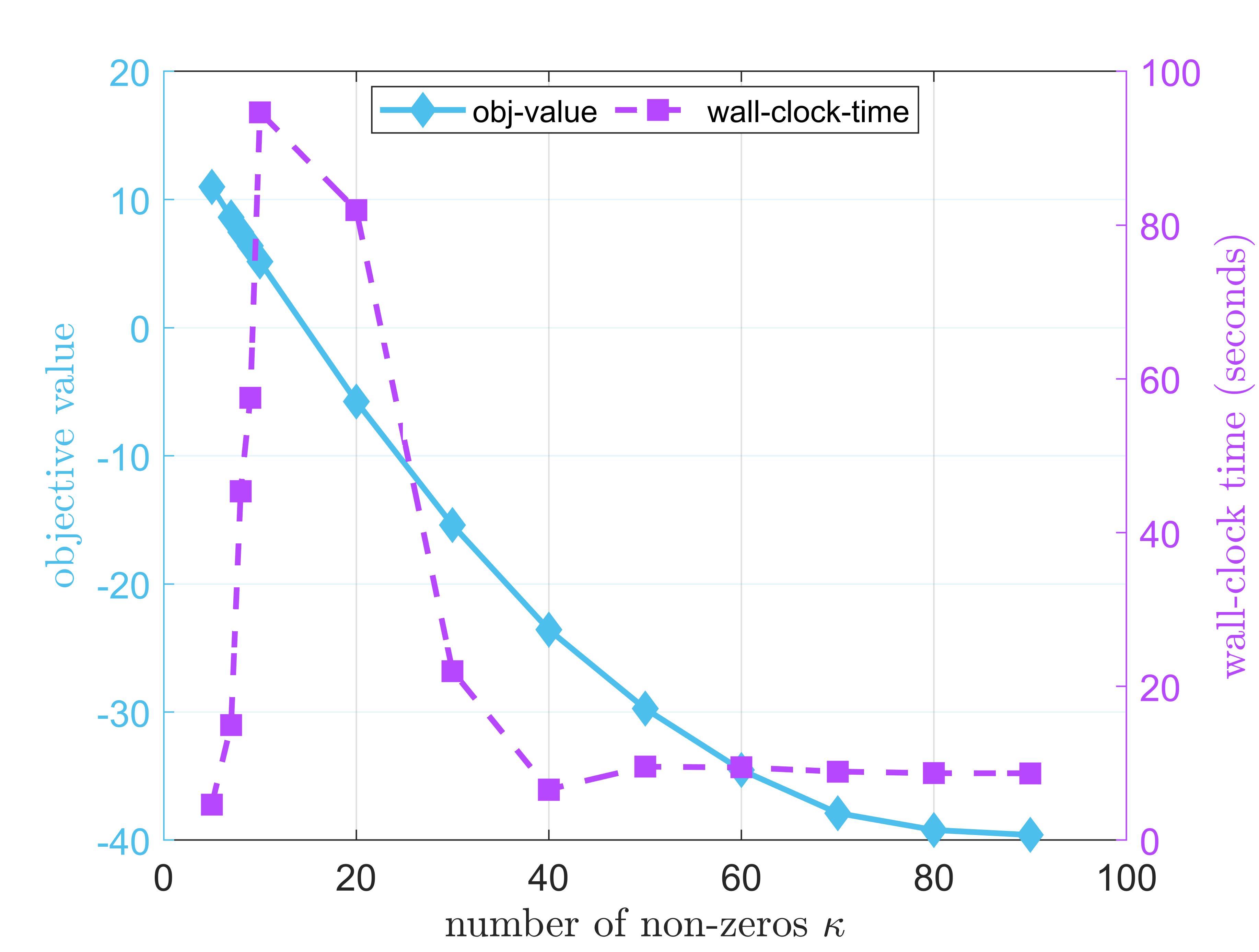}
    \caption{\footnotesize Objective value and wall-clock time for different values of $\kappa$ for $n = 100$; the left axis shows the behavior of the optimal objective value and the right axis represents the execution time.}  
    \label{sc1: qc}
\end{figure}

\begin{figure}
    \centering
    \includegraphics[width=0.55\textwidth]{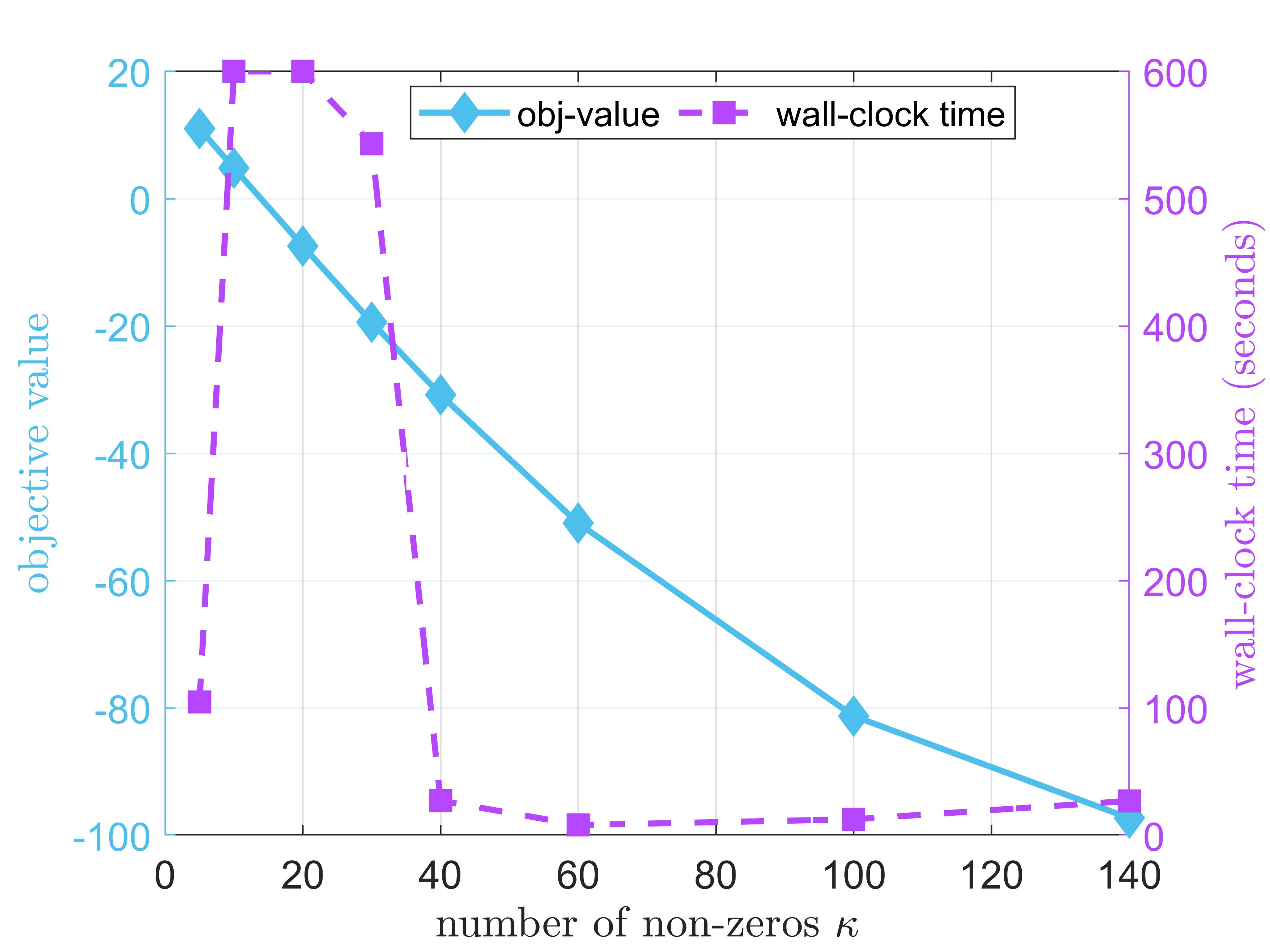}
    \caption{\footnotesize Objective value and wall-clock time for different values of $\kappa$ for $n = 200$; the left axis shows the behavior of the optimal objective value and the right axis represents the execution time.}
    \label{sc2: qc}
\end{figure}

In this experiment, two main scenarios with different problem sizes are considered. The first scenario has $n = 100$ variables while the second has  $n = 200$. In both scenarios, we evaluate DiPOA according to a different number of nonzero elements $\kappa$.
 \begin{table}[]
    \centering
    \caption{\footnotesize Numerical Results for SQCP with $n = 100$ Variables \label{sc: qc1}}
    {    \begin{tabular}{SSSSS}
    \toprule
            {$\kappa$} &{ \Verb!objective-value!} & { \Verb!rel-gap[$\%$]!} & { \Verb!total-cuts!} & { \Verb!time[sec]!}\\\toprule
            {$5 $} &  {$10.994 $} &  {$ 0.131$} &  {$30 $} &  {$4.64 $} \\
            {$7 $} &  {$8.601 $} &  {$0.143 $} &  {$ 80$} &  {$ 14.96$} \\
            {$8 $} &  {$7.457 $} &  {$ 0.142$} &  {$ 180$} &  {$ 45.388$} \\
            {$9 $} &  {$6.378 $} &  {$ 0.146$} &  {$ 200$} &  {$57.519 $} \\
            {$10 $} &  {$5.141 $} &  {$ 0.146$} &  {$270 $} &  {$94.68 $} \\
            {$20 $} &  {$-5.759 $} &  {$0.139 $} &  {$ 250$} &  {$81.962 $} \\
            {$30 $} &  {$-15.411 $} &  {$ 0.144$} &  {$ 90$} &  {$21.986 $} \\
            {$40 $} &  {$-23.576 $} &  {$0.127 $} &  {$ 50$} &  {$6.604 $} \\
            {$50 $} &  {$-29.733 $} &  {$0.144 $} &  {$ 50$} &  {$ 9.554$} \\
            {$60 $} &  {$-34.529 $} &  {$0.101 $} &  {$ 40$} &  {$9.462 $} \\
            {$70 $} &  {$-37.885 $} &  {$ 0.122$} &  {$ 30$} &  {$8.921 $} \\
            {$80 $} &  {$-39.208 $} &  {$ 0.060$} &  {$ 30$} &  {$8.712 $} \\
            {$90 $} &  {$-39.595 $} &  {$ 0.054$} &  {$ 30$} &  {$8.699 $} \\
    \bottomrule
    \end{tabular}
    }
    {}
\end{table}
The numerical results of the first scenario appear in Table \ref{sc: qc1}. It can be noticed that, by decreasing the degree of the solution sparsity, the optimal objective value also decreases. This behavior is expected since, by reducing the sparsity of the solution, more variables can appear in the solution and therefore the optimal objective approaches the value of the dens problem (\textit{i.e.}, problem \eqref{dis-ccp} without sparsity constraint). 
The number of cuts and hence the execution time also changes depending on the sparsity of the solution. In particular, as the number of nonzero variables $\kappa$ increases, more iterations are needed by the DiPOA algorithm to converge. However, this phenomenon occurs until $\kappa = \kappa_{\max}$, where $\kappa_{\max}$ is the largest number of nonzero variables for which the highest number of iterations is needed. For $\kappa > \kappa_{\max}$ the complexity of the problem decreases. The behavior of the objective value and also the solution time is illustrated in Figure \ref{sc1: qc}, according to the sparsity of the solution. For this problem instance, the figure shows that $\kappa_{\max} = 10$, at which sparsity-level the wall-clock time is $94.68$~s and the total number of cuts is $270$. 
Finally, the numerical results of the second scenario are provided in Table \ref{sc: qc2}. As mentioned before, in this case, we scale the problem and set $n = 200$.  
 \begin{table}[]
    \centering
    \caption{\footnotesize Numerical Results for SQCQP with $n = 200$ \label{sc: qc2}}
    {    \begin{tabular}{SSSSS}
    \toprule
            {$\kappa$} &{ \Verb!objective-value!} & { \Verb!rel-gap[$\%$]!} & { \Verb!total-cuts!} & { \Verb!time[sec]!}\\\toprule
            {$5 $} &  {$ 11.001$} &  {$0.14 $} &  {$240 $} &  {$ 104.768$} \\
            {$ 10$} &  {$4.816 $} &  {$1.34 $} &  {$ 810$} &  {$ 600$} \\
            {$ 20$} &  {$ -7.427$} &  {$1.27 $} &  {$880 $} &  {$ 600$} \\
            {$30 $} &  {$ -19.397$} &  {$ 0.11$} &  {$740 $} &  {$543.268 $} \\
            {$40 $} &  {$ -30.802$} &  {$ 0.14$} &  {$ 80$} &  {$27.133 $} \\
            {$60 $} &  {$ -50.985$} &  {$ 0.14$} &  {$20 $} &  {$ 8.021$} \\
            {$100 $} &  {$-81.291 $} &  {$ 0.12$} &  {$40 $} &  {$12.498 $} \\
            {$140 $} &  {$ -97.308$} &  {$ 0.05$} &  {$ 50$} &  {$26.952$} \\
           
    \bottomrule
    \end{tabular}
    }
    {}
\end{table}
The behavior of the objective value and the wall-clock time is illustrated in Figure \ref{sc2: qc}. As expected, a similar behavior emerges from the application of DiPOA, however, a large computational time is needed in some instances. 

We also compared DiPOA against Gurobi for solving the convex
SQCQP. In most problem instances, Gurobi \cite{gurobi} achieved better performance for the centralized case in terms of wall clock time, which is expected considering the communication and computing overhead of DiPOA. Regarding the distributed setup, Gurobi like other centralized solvers cannot be directly applied because the objective-defining data is private to each node and spread over the CN.
\subsection{Discussion}
In this section, we evaluated the performance of DiPOA for DSLR and SQCQP problems. We considered $70$ problem instances under different scenarios and settings. According to the numerical experiments, DiPOA performance and efficiency were acceptable for both problems in most situations. The numerical results and comparison with state-of-the-art MINLP solvers also showed that, for the conditions when distributed computations are inevitable, DiPOA provides a reliable and optimal solution. 

\textcolor{black}{As a general conclusion, the separability of the optimization problem can be used by distributed algorithms to speed up the computation by splitting the computational burden among several processors. Essentially, each processor would be responsible to solve a portion of the problem independent of other processors and in parallel. For small-scale problems, however, a distributed algorithm might not necessarily outperform centralized algorithms taking into account process synchronization and inter-process communication. 
 A distributed algorithm can be a suitable choice in scenarios where either the problem is inherently distributed or large-scale. The former case refers to the problems where the data is spread over a large network and data privacy matters. In this case, because of the data privacy and the size of the network, it is not usually practical to process the problem data in a central node. In the latter case, due to hardware limitations, it is difficult (if not impossible) to process the entire problem data in a single computational node and a distributed algorithm over multiple nodes can be useful.}

\section{Conclusions}\label{concolusions}
This work proposed an algorithmic methodology to solve SCP problems consisting of convex and differentiable nonlinear functions. Employing the idea of LFCs and MINLP modeling framework, the SCP problem was reformulated as a distributed consensus MINLP optimization problem.  To solve the consensus MINLP problem, we developed a distributed version of the OA algorithm whereby the primal subproblems are solved in a fully decentralized manner. Moreover, we introduced a second-order optimality cut generation methodology, along with a specialized feasibility pump heuristic, to improve the convergence speed of the proposed DiPOA algorithm. Finally, the performance of the DiPOA algorithm was evaluated by solving several instances of  DSLR and SQCQP problems under different scenarios.

As for future work, we intend to improve DiPOA performance by implementing single-tree optimization and designing fast primal heuristics. Another direction for research is the development of a decentralized solution for the master's problem, possibly utilizing Lagrangean and stochastic decomposition strategies.



%
\section*{Funding}

This work was funded in part by  Funda\c{c}\~ao de Amparo \`a Pesquisa e Inova\c{c}\~ao do Estado de Santa Catarina  (FAPESC) under grant 2021TR2265 and Coordena{\c{c}}{\~a}o de Aperfei{\c{c}}oamento de Pessoal
de  N{\'i}vel Superior (CAPES, Brazil) under the project PrInt CAPES-UFSC 698503P1.

\section*{Statements and Declarations}

\textbf{Competing Interests:} The authors declare that they have no conflict of interest.

\bibliographystyle{spmpsci}      
\bibliography{Refs.bib}   


\end{document}